\documentclass[11pt,a4paper]{article}
\usepackage{latexsym}
\usepackage{amsfonts}
\usepackage{amssymb}

\begin{document}

\setlength{\belowdisplayshortskip}{1\belowdisplayskip}
\newlength{\indentation} \setlength{\indentation}{-1\parindent}
\newcommand{\no}{\hspace*{\indentation}}
\renewcommand{\arraystretch}{1.75}
\setlength{\arrayrulewidth}{0.5\arrayrulewidth}

\newcounter{propcount}
\newenvironment{proplist}[2]{\begin{list}{\normalfont{(\roman{propcount})}}
    {\usecounter{propcount}
    \setlength{\topsep}{4pt} \setlength{\itemsep}{#2pt}
    \setlength{\parsep}{2pt} \setlength{\labelsep}{3mm}
    \setlength{\leftmargin}{#1mm}
    \setlength{\rightmargin}{#1mm}}}
       {\end{list}}

\newcounter{prop2count}
\newenvironment{prop2list}[3]{\begin{list}{\normalfont{(\roman{prop2count})}}
    {\usecounter{prop2count}
    \setlength{\topsep}{4pt} \setlength{\itemsep}{#3pt}
    \setlength{\parsep}{2pt} \setlength{\labelsep}{3mm}
    \setlength{\leftmargin}{#1mm}
    \setlength{\rightmargin}{#2mm}}}
       {\end{list}}

\newlength{\tablength} \setlength{\tablength}{2.5mm}
\newcommand{\tab}[1]{\hspace*{#1\tablength}}
\newcommand{\nm}[1]{\mbox{\ensuremath{\| #1 \|}}}
\newcommand{\implies}{\ensuremath{\Rightarrow}}
\newcommand{\fa}{\ensuremath{\ \, \forall \,}}
\renewcommand{\iff}{\ensuremath{\; \Longleftrightarrow \;}}
\renewcommand{\sp}[1]{\ensuremath{\mathrm{sp}(#1)}}
\newcommand{\gap}[1]{\ensuremath{\, #1 \,}}
\newcommand{\inset}[2]{\ensuremath{\{ #1 \, | \, #2 \} }}
\newcommand{\inprod}[2]{\ensuremath{\langle #1 , #2 \rangle}}
\newcommand{\anti}[2]{\ensuremath{ \{ #1 , #2 \} }}
\newcommand{\CC}{\ensuremath{\mathbb{C}}}
\newcommand{\RR}{\ensuremath{\mathbb{R}}}
\newcommand{\TT}{\ensuremath{\mathbb{T}}}
\newcommand{\SB}{\ensuremath{\mathbb{S}}}
\newcommand{\ZZ}{\ensuremath{\mathbb{Z}}}
\newcommand{\HH}{\ensuremath{\mathbb{H}}}
\newcommand{\NN}{\ensuremath{\mathbb{N}}}
\newcommand{\iny}{\ensuremath{\infty}}
\newcommand{\bd}[1]{ \ensuremath{r_{\NN}(\ensuremath{#1})}}
\newcommand{\bp}{\ensuremath{p}}

\newcommand{\noline}{\vspace*{1\parskip}}
\newcommand{\preskip}{\vspace*{1\belowdisplayskip}}
\newcommand{\display}[1]{$$#1$$}
\newcommand{\smdisplay}[1]{\\[4pt] \no \centerline{#1} \\[4pt]}
\newcommand{\text}[1]{\mbox{\normalfont #1}}
\newcommand{\mod}[1]{\ensuremath{
        \text{ \hspace*{-3pt} mod \hspace*{-3pt} } #1}}

\newtheorem{dfnS}{Definition}[section]
\newtheorem{dfn-lemmaS}[dfnS]{Definition-Lemma}
\newtheorem{dfn-corS}[dfnS]{Definition-Corollary}
\newtheorem{lemmaS}[dfnS]{Lemma}
\newtheorem{propS}[dfnS]{Proposition}
\newtheorem{thmS}[dfnS]{Theorem}
\newtheorem{corS}[dfnS]{Corollary}
\newtheorem{exampleS}[dfnS]{Example}
\newtheorem{remarkS}[dfnS]{Remark}
\newcommand{\EXS}[1]{\begin{exampleS} {\normalfont #1} \end{exampleS}}
\newcommand{\REMS}[1]{\begin{remarkS} {\normalfont #1} \end{remarkS}}

\newtheorem{dfn}{Definition}[subsection]
\newtheorem{dfn-cor}[dfn]{Definition-corollary}
\newtheorem{lemma}[dfn]{Lemma}
\newtheorem{prop}[dfn]{Proposition}
\newtheorem{thm}[dfn]{Theorem}
\newtheorem{cor}[dfn]{Corollary}
\newtheorem{example}[dfn]{Example}
\newtheorem{remark}[dfn]{Remark}
\newcommand{\EX}[1]{\begin{example} {\normalfont #1} \end{example}}
\newcommand{\REM}[1]{\begin{remark} {\normalfont #1} \end{remark}}
\newcommand{\proof}[1]{\normalfont{\noindent \textbf{Proof} \  \  #1 \hfill $\Box$} \par}

        \title{Dixmier Traces as Singular Symmetric Functionals
and Applications to Measurable Operators}
\author{Steven Lord\footnote{Research supported by the Australian Research
Council}\ , Aleksandr Sedaev\footnote{Research supported by RFFI
grant 02-01-00146 and by the
 Scientific program \lq\lq Universities of Russia \rq\rq grant UR 04.01.051}
   \  and   Fyodor Sukochev
}
\date{}

                \maketitle

\begin{abstract}
\noindent We unify various constructions and contribute to the
theory of singular symmetric functionals on Marcinkiewicz
function/operator spaces. This affords a new approach to the
non-normal Dixmier and Connes-Dixmier traces (introduced by
Dixmier and adapted to non-commutative geometry by Connes) living
on a general Marcinkiewicz space associated with an arbitrary
semifinite von Neumann algebra. The corollaries to our approach,
stated in terms of the operator ideal $\mathcal{L}^{(1,\iny)}$
(which is a special example of an operator Marcinkiewicz space),
are: (i) a new characterization of the set of all positive
measurable operators from $\mathcal{L}^{(1,\iny)}$, i.e. those on
which an arbitrary Connes-Dixmier trace yields the same value. In
the special case, when the operator ideal $\mathcal{L}^{(1,\iny)}$
is considered on a type $I$ infinite factor, a bounded operator
$x$ belongs to $\mathcal{L}^{(1,\iny)}$ if and only if the
sequence of singular numbers $\{s_n(x)\}_{n \ge 1}$ (in the
descending order and counting the multiplicities) satisfies
$\|x\|_{(1,\infty)}:= \sup_{N \geq 1}
\frac{1}{Log(1+N)}\sum_{n=1}^Ns_n(x)<\infty$. In this case, our
characterization amounts to saying that a positive element $x\in
\mathcal{L}^{(1,\iny)}$ is measurable if and only if\linebreak
$\lim_{N\to\infty} \frac{1}{LogN}\sum_{n=1}^Ns_n(x)$ exists; (ii)
the set of Dixmier traces and the set of Connes-Dixmier traces are
norming sets (up to equivalence) for the space
$\mathcal{L}^{(1,\iny)}/\mathcal{L}^{(1,\iny)}_0$, where the space
$\mathcal{L}^{(1,\iny)}_0$ is the closure of all finite rank
operators in $\mathcal{L}^{(1,\iny)}$ in the norm $\| . \|_{(1,\infty)}$.
\end{abstract}

\section*{Introduction}

In \cite{Dix} Dixmier proved the existence of non-normal traces on
the von Neumann algebra $B(H)$. Dixmier's original construction
involves singular dilation invariant positive linear functionals
$\omega$ on $\ell^\iny(\NN)$. This construction was altered by A.
Connes \cite{CN} (see also Definition 5.2 below) who defined
non-normal traces via the composition of the Cesaro mean and a
state on $C_b([0,\iny)) / C_0([0,\iny))$. In \cite{DPSS},
\cite{DPSSS} and \cite{DPSSS2} the traces of Dixmier in \cite{Dix}
were broadly generalized as singular symmetric functionals on
Marcinkiewicz function (respectively, operator) spaces $M(\psi)$
on $[0,\iny)$ (respectively, on a semifinite von Neumann algebra).
The symmetric functionals in \cite{DPSSS} and \cite{DPSSS2}
involve Banach limits, that is, singular translation invariant
positive linear functionals $L'$ on $\ell^\iny(\NN)$. We extend
the construction of Dixmier in Definition 1.7 and Connes
in Definition 5.2 (verified in Theorem 6.3) by extending
the notion of Banach limits to $C_b([0,\iny))$.

The identification of the commutative specialization of \mbox{\text{(Connes-)}}Dixmier traces
as singular symmetric
functionals has some pivotal consequences.
The established
theory of Banach limits \cite{GL} and singular symmetric
functionals on Marcinkiewicz spaces \cite{DPSS}, \cite{DPSSS},
\cite{DPSSS2} can be applied to questions concerning the
\mbox{\text{(Connes-)}} Dixmier trace, a central notion in Connes'
non-commutative geometry \cite{CN}. Conversely, ideas in Connes'
non-commutative geometry, such as measurability of operators
\cite[IV.2.$\beta$, Definition 7]{CN}, lend themselves to
generalization to abstract Marcinkiewicz spaces (Definition 3.2
and Definition 3.5). As a result, we have been able to present a
new characterization of measurable operators (see Theorem 5.12,
Remark 5.13 and Theorem 6.6).

\medskip \noindent The paper is structured as follows.

Section 1 introduces
Banach limits, almost convergence (extending the notions of G. Lorentz \cite{GL})
and the theory of singular symmetric functionals on
the Marcinkiewicz space $M(\psi)$ defined by a concave function $\psi$
\cite{DPSS}, \cite{DPSSS}.
The construction of singular symmetric functionals on $M(\psi)$
\cite{DPSSS} (Definition 1.6 below) is extended by Definition 1.7.

Section 2 introduces sufficient conditions to identify
the singular symmetric functionals of \cite{DPSSS}
with those of Definition 1.7, see Theorem 2.3 and Theorem 2.7.
A result in \cite{DPSSS}, on
the Riesz semi-norm of a function $x$ in a Marcinkiewicz space $M(\psi)$
as the supremum of the values $\{ f(x) \}$ where $\{ f \}$ is a set of
singular symmetric functionals on $M(\psi)$, is extended in Theorem 2.8.

Section 3 contains
an analysis of various notions of a measurable element of a
Marcinkiewicz space $M(\psi)$, introduced in Definitions 3.2 and
3.5, and their coincidence (Theorem 3.7 and Corollary 3.9, see
also Theorem 3.14).

The results of Section 2 and Section 3 concern singular
symmetric functionals on $M(\psi)$ parameterised by the set
of strictly increasing, invertible, differentiable
and unbounded functions $\kappa: [0,\infty) \to [0,\infty)$.
Section 4 summarises the conditions on the function
$\kappa$ used in Section 2 and Section 3.  Theorem 4.4,
which extends the existence results of \cite{DPSS},
demonstrates an equivalence between the growth
of the concave function $\psi$ and the existence
of functions $\kappa$ which satisfy the hypotheses
of results in Section 2 and 3.

A subset of the collection of extended Banach
limits, called Cesaro-Banach limits (Definition 5.4) is studied
further in Section 5.  It is demonstrated that this subset is coincident
with
the generalized limits employed by Connes to construct the
\mbox{Connes-Dixmier} traces used in non-commutative geometry.
Theorem 5.6 identifies (the commutative specialization of)
\mbox{Connes-Dixmier} traces as a sub-class of the singular
symmetric functionals studied in \cite{DPSSS}, \cite{DPSSS2}.
Results on \mbox{Connes-Dixmier} traces then follow from the
general theory of singular symmetric functionals on Marcinkiewicz
spaces developed in the preceding sections (Theorem 5.12).

Section 6 considers the special example of the Marcinkiewicz
space $M(\psi)$ where $\psi(t)=\log(1+t)$
(recognized from non-commutative geometry
as the space $\mathcal{L}^{(1,\infty)}$). Here, we
summarize and present our results (Theorems 6.1, 6.2, 6.3, 6.4 and 6.6)
for Dixmier and Connes-Dixmier traces on the operator
Marcinkiewicz spaces associated with semifinite von Neumann
algebras of type \emph{I} and \emph{II}.  In particular, Theorems 6.1, 6.2, 6.3, 6.4
and 6.6 apply to the operator ideals and traces
of non-commutative geometry \cite{CN}.

\section{Preliminaries}

\subsection{Banach Limits, Almost Convergence and Almost
Piecewise Linearity}

\medskip \noindent Let $H$ be one of the semigroups $\NN := \{1,2,... \}$ or
$\RR_+ := [0,\iny)$ equipped with the topology and order induced
by the locally compact additive group $\RR$. Let $C_b(H)$ be the
space of bounded continuous functions on $H$. Define the
translation operator \display{T_s(f)(t) = f(t+s) \ \fa s,t \in H,
f \in C_b(H).} An element $L \in C_b(H)^*$ is called translation
invariant if \display{L(T_s(f)) = L(f) \fa s \in H, f \in C_b(H).}
A \textbf{Banach limit} $L$ on $C_b(H)$ is a translation invariant
positive linear functional on $C_b(H)$ such that $L(1) = 1$. This
extends the notion of a Banach limit investigated in \cite{GL} in
the context of the semigroup $\NN$ of all natural numbers. Let
$BL(H)$ denote \textbf{the set of all Banach limits on} $C_b(H)$.
It is easy to see that every $L \in BL(H)$ vanishes on compactly
supported elements from $C_b(H)$ and that
$$
\liminf_{t \to \iny} f(t) \leq L(f) \leq \limsup_{t \to \iny} f(t)
$$
for any positive $f \in C_b(H)$.

We extend the notion of almost convergent sequences \cite{GL}.

\begin{dfnS}
A function $f \in C_b(H)$ is said to be \textbf{\emph{almost convergent
at infinity}} if $L_1(f) = L_2(f) \fa L_1,L_2 \in BL(H)$.
\end{dfnS}

Let $f \in C_b(H)$ be almost convergent at infinity.  We denote the value $A :=
L(f) \fa L \in BL(H)$ by
$$
\text{F-}\lim f = A
$$
following
G. Lorentz \cite{GL}.  In particular we write
 ${\text{F-}\lim_{n\to \infty}
a_n}$ for $\alpha = \{a_n \}_{n=1}^\iny \in \ell^\iny(\NN)$ and
${\text{F-}\lim_{t\to \infty}
g(t)}$ for $g \in C_b([0,\iny))$.

\bigskip \noindent
Let $\alpha = \{a_n \}_{n=1}^\iny \in \ell^\iny(\NN)$. Let
$\chi_E$ be the characteristic function for $E \subset [0,\iny)$.
Define the \textbf{piecewise linear extension map} \display{\bp :
\ell^\iny(\NN) \to C_b([0,\iny))} by \display{\bp(\alpha)(t) =
\sum_{n=0}^\infty  \Big( a_n + (a_{n+1}-a_n)(t-n) \Big)
\chi_{[n,n+1)}(t),}  where $a_0=0$ by definition. The following
lemma is an elementary application of the definition, hence the
proof is omitted.

\begin{lemmaS}
The map $\bp : \ell^\iny(\NN) \to C_b([0,\iny))$ is a positive
linear isometry with the following properties
\begin{prop2list}{10}{2}{3}
\item $\bp(1_{\ell^\iny}) = 1$,
\item $\nm{\bp(\alpha)} = \nm{\alpha}$ for all $\alpha \in \ell^\iny(\NN)$,
\item $T_{k}(\bp(\alpha)) = \bp(T_k(\alpha))$
for all $\alpha \in \ell^\iny(\NN)$ and $k \in \NN$.
\end{prop2list}
\end{lemmaS}

\noindent Let $g \in C_b([0,\iny))$. Define the
\textbf{restriction map} $r_{\NN}$ and \textbf{averaging map}
$E_{\NN}$, acting from $ C_b([0,\iny))$ onto $\ell^\iny(\NN)$ by
\display{\bd{g} := \{ g(n) \}_{n=1}^\iny \ \ , \ \ E_\NN(g) := \{
\int_{n-1}^{n} g(s)ds \}_{n=1}^\iny.} The following lemma is an
elementary application of the definitions.

\begin{lemmaS}
The maps $r_{\NN}, E_{\NN} : C_b([0,\iny)) \to \ell^\iny(\NN)$ are positive
linear surjections with the following properties
\begin{prop2list}{10}{2}{3}
\item $\bd{1} = E_{\NN}(1) = 1_{\ell^\iny}$,
\item $\nm{\bd{g}} \leq \nm{g}$ and $\nm{E_{\NN}(g)} \leq \nm{g}$
for all $g \in C_b([0,\iny))$,
\item $\bd{T_{a + k}(g)} = T_k \bd{T_a(g)}$ and
$E_{\NN}(T_{a + k}(g)) = T_k E_{\NN}(T_a(g))$ \\
\hspace*{0.4cm} for all $a \in [0,\iny)$, $g \in C_b([0,\iny))$ and $k \in \NN$.
\end{prop2list}
\end{lemmaS}

\noindent The following notion shall become an important concept in Section 2.

\begin{dfnS}
Let $g \in C_b([0,\iny))$.  We say $g$ is \textbf{\emph{almost piecewise linear
at infinity}} if $L(g - \bp \bd{g}) = 0 \fa L \in BL(\RR_+)$.
\end{dfnS}

\subsection{Singular symmetric functionals on Marcinkiewicz
spaces}

\noindent We introduce the notation of \cite{DPSSS}. Let $m$ be
the Lebesgue measure on $[0,\iny)$. Let $x$ be a measurable
function on $[0,\iny)$. Define the \textbf{decreasing
rearrangement of} $x$ by \display{x^*(t) = \inf \inset{s \geq
0}{m(\{ |x| > s\}) \leq t} , \ t > 0.} Let $\Omega_\iny$ denote
the set of concave functions $\psi : [0,\iny) \to [0,\iny)$ such
that $\lim_{t \to 0^+} \psi(t) = 0$ and $\lim_{t \to \iny} \psi(t)
= \iny$. Important functions belonging to $\Omega_\iny$ include
$t$, $\log(1+t)$, $t^\alpha$ and $(\log(1+t))^{\alpha}$ for $0 <
\alpha < 1$, and $\log(1+\log(1+t))$. Let $\psi \in \Omega_\iny$.
Define the \textbf{weighted mean function} \display{\phi(x)(t) :=
\frac{1}{\psi(t)} \int_0^{t} x^*(s) ds,\ t>0} and denote by
$M(\psi)$ the \textbf{Marcinkiewicz space} of measurable functions
$x$ on $[0,\iny)$ such that \display{\nm{x}_{M(\psi)} := \sup_{t
> 0} \phi(x)(t) = \nm{\phi(x)}_\iny
 < \iny.}

\noindent The norm closure of $M(\psi) \cap L^1([0,\iny))$ in
$M(\psi)$ is denoted by $M_1(\psi)$. For every $\psi\in
\Omega_\iny$, we have $M_1(\psi)\neq M(\psi)$. We define the
\textbf{Riesz semi-norm} on $M(\psi)$ by  \display{\rho_1(x) :=
\inf \inset{ \nm{x-y}_{M(\psi)}}{y \in M_1(\psi)}=\limsup _{t\to
\infty}\phi(x)(t),} (see \cite[Proposition 2.1]{DPSSS}). The
Banach space $(M(\psi),\nm{.}_{M(\psi)})$ is an example of a
rearrangement invariant space \cite{LT}, also termed a symmetric
space \cite{KPS}. Let $M_+(\psi)$ denote the set of positive
functions of $M(\psi)$.

\begin{dfnS}

A positive homogeneous functional $f:M_+(\psi)\to [0,\infty)$ is
\text{(i)} \textbf{\emph{symmetric}} if $f(x) \leq f(y)$ for all $x,y \in M_+(\psi)$ such
that $\int_0^t x^*(s)ds \leq \int_0^t y^*(s)ds \fa t \in
[0,\iny)$, and \text{(ii)} \textbf{\emph{supported at infinity}}, or \textbf{\emph{singular}} on
$M(\psi)$, if $f(|x|) = 0$ for all $x \in M_1(\psi)$.

\end{dfnS}
If such a functional is {\it additive}, then it
can be extended by linearity to a bounded linear positive
functional on $M(\psi)$. Let $M_+(\psi)_{\mathrm{sym},\iny}^*$
denote the \textbf{cone of additive symmetric functionals on
$M_+(\psi)$ supported at infinity}, \cite[Section 2]{DPSSS}. Not
every Marcinkiewicz space $M(\psi)$, $\psi\in\Omega_\iny$, admits
non-trivial additive singular symmetric functionals. %
Necessary
and sufficient conditions for the existence of such functionals on
a Marcinkiewicz space $M(\psi)$ may be found in \cite[Theorem
3.4]{DPSS} and will be considered in Theorem 2.7 and Theorem 4.4
below.

\bigskip \noindent
Let $\kappa : [0,\iny) \to [0,\iny)$ be an increasing, continuous
and unbounded function.  Define the $\kappa$\textbf{-weighted mean
function on} $M(\psi)$
$$
\phi_\kappa(x)(t) := \phi(x)(\kappa(t)) =
\frac{1}{\psi(\kappa(t))} \int_0^{\kappa(t)} x^*(s) ds,\quad t>0.
$$
Then, as $\phi(x) \in C_b([0,\iny))$ for each $x \in M(\psi)$, we
have $\phi_\kappa(x) \in C_b([0,\iny))$ for each $x \in M(\psi)$
and the sequences \display{ \bd{\phi(x)} = \{ \phi(x)(n)
\}_{n=1}^\iny \ \ \ \text{and} \ \ \ \bd{\phi_\kappa(x)} = \{
\phi_\kappa(x)(n)) \}_{n=1}^\iny} are bounded.

\begin{dfnS}
Let $\psi \in \Omega_\iny$ and $\kappa : [0,\iny) \to [0,\iny)$ be
increasing, continuous and unbounded. Let $x \in M_+(\psi)$ and
$L' \in BL(\NN)$. Define \display{f_{L',\kappa}(x) := L' (
\bd{\phi_\kappa(x)}) = L' \Big( \{ \phi(x)(\kappa(n))
\}_{n=1}^\iny \Big) .}
\end{dfnS}

\noindent In \cite{DPSSS} necessary and sufficient conditions were
found on the sequence $\{ \kappa(n) \}_{n=1}^\iny$ and the
function $\psi\in \Omega_\iny$ such that $f_{L',\kappa} \in M_+(\psi)_{\mathrm{sym},\iny}^*$
for all $L' \in BL(\NN)$.
It is natural to introduce the following extension.

\begin{dfnS}
Let $\psi \in \Omega_\iny$ and $\kappa : [0,\iny) \to [0,\iny)$ be
increasing, continuous and unbounded. Let $x \in M_+(\psi)$ and $L
\in BL(\RR_+)$.  Define \display{f_{L,\kappa}(x) :=
L(\phi_\kappa(x)). }
\end{dfnS}

The analysis of the functionals $f_{L,\kappa}$ on $M_+(\psi)$ begins in Section 2.2.
We finish the preliminaries with the following proposition and remark.

\begin{propS}
\ Let $\psi \in \Omega_\iny$, $L \in BL(\RR_+)$ and $ L' \in
BL(\NN)$. Let $\kappa_1, \kappa_2 : [0,\iny) \to [0,\iny)$ be increasing continuous
and unbounded functions such that
$\kappa_1-\kappa_2$ is bounded.  Then
$f_{L,\kappa_1}(x)=f_{L,\kappa_2}(x)$ and
$f_{L',\kappa_1}(x)=f_{L',\kappa_2}(x)$ for all $x\in
M_+(\psi)$.
\par
\preskip
\proof{Since
$\nm{x}_{M(\psi)} \geq \psi(t)^{-1}\int_0^tx^*(s)ds \ge x^*(t)t\psi(t)^{-1}$
for all $t>0$, we have
$$
-\|x\|_{M(\psi)}\frac{\psi'(t)}{\psi(t)}  \le
-\frac{\int_0^tx^*(s)ds\psi'(t)}{\psi^2(t)} \le \phi(x)'(t) \le
\frac{x^*(t)}{\psi(t)} \le \frac{\|x\|_{M(\psi)}}{t}.
$$
Hence for any $t > 0$,
$$
- \nm{x}_{M(\psi)} \log(\psi(t'))|_{\kappa_1(t)}^{\kappa_2(t)} \leq
    \phi(x)(t') |_{\kappa_1(t)}^{\kappa_2(t)} \leq \nm{x}_{M(\psi)}
 \log(t')|_{\kappa_1(t)}^{\kappa_2(t)}
$$
or
$$
-\nm{x}_{M(\psi)} \log \frac{\psi(\kappa_2(t))}{\psi(\kappa_1(t))} \leq
\phi_{\kappa_2}(x)(t) - \phi_{\kappa_1}(x)(t) \leq \nm{x}_{M(\psi)}
\log \frac{\kappa_2(t)}{\kappa_1(t)} \eqno{(1.1)}
$$
Let $f$ be an unbounded concave function.  Then
$\Big|\frac{f(\kappa_2(t))}{f(\kappa_1(t))} - 1 \Big| =
\frac{|f(\kappa_2(t)) - f(\kappa_1(t))|}{f(\kappa_1(t))} \leq A
\frac{|\kappa_2(t) - \kappa_1(t)|}{f(\kappa_1(t))} \leq \frac{A
B}{f(\kappa_1(t))}$ for $A,B > 0$ and $t$ sufficiently large by
the hypothesis $f$ is concave and $\kappa_2 - \kappa_1$ is
bounded. Hence $\lim_{t \to \iny}
\frac{f(\kappa_2(t))}{f(\kappa_1(t))} = 1$. Then $\lim_{t \to
\iny} \log \Big| \frac{\kappa_2(t)}{\kappa_1(t)} \Big| = \lim_{t
\to \iny} \log \Big| \frac{\psi(\kappa_2(t))}{\psi(\kappa_1(t))}
\Big| = 0$ and $\phi_{\kappa_2}(x)(t) - \phi_{\kappa_1}(x)(t) \in
C_0([0,\iny))$ by (1.1). Since $L\in BL({\Bbb R}_+)$
(respectively, $ L'\in BL({\Bbb N}))$ vanishes on functions
(respectively, sequences) tending to 0 at infinity, we conclude
that $f_{L,\kappa_2}(x) - f_{L,\kappa_1}(x) = L(\phi_{\kappa_2}(x)
- \phi_{\kappa_1}(x)) = 0$ (respectively, $f_{L',\kappa_2}(x)
-f_{L',\kappa_1}(x) = L'(\{\phi_{\kappa_2}(n) -
\phi_{\kappa_1}(n)\}) = 0$). }
\end{propS}

\REMS{Proposition 1.8 introduces the notion of
equivalence classes of continuous increasing unbounded functions
that result in the same functional on $M_+(\psi)$.

Let $\psi \in \Omega_\iny$ and $\kappa_1, \kappa_2 : [0,\iny) \to
[0,\iny)$ be continuous increasing and unbounded functions. We
define an equivalence relation $\sim_\psi$ by \display{\kappa_1
\sim_\psi \kappa_2 \text{ \ if \ } f_{L,\kappa_1}(x) =
f_{L,\kappa_2}(x) \fa L \in BL(\RR_+) \fa x \in M_+(\psi).} Let
$[\kappa]$ denote the equivalence class, with respect to the
relation $\sim_\psi$, of a continuous increasing and unbounded
function $\kappa : [0,\iny) \to [0,\iny)$. It easily follows from
Proposition 1.8 that the class $[\kappa]$ contains a strictly
increasing, invertible, unbounded function $\hat{\kappa}$ such
that $\hat{\kappa}(0)=0$.  The function $\hat{\kappa}$ can be
chosen to be differentiable or even piecewise linear if required.
Hence, to analyse all functionals $f_{L,\kappa}$ on $M_+(\psi)$
where $\kappa : [0,\iny) \to [0,\iny)$ is a continuous increasing
unbounded function, it is sufficient to consider the set
$\mathcal{K}$ \textbf{of strictly increasing, invertible,
differentiable, unbounded functions} $\hat{\kappa} : [0,\iny) \to
[0,\iny)$ \textbf{such that} $\hat{\kappa}(0)=0$.

}

\section{Symmetric Functionals involving Banach Limits}

This section demonstrates that: (i) the sets of functionals
$\inset{f_{L',\kappa}}{L' \in BL(\NN)}$ (Definition 1.6) and
$\inset{f_{L,\kappa}}{L \in BL(\RR_+)}$ (Definition 1.7) provide
the same set of functionals on $M_+(\psi)$ supported at infinity
for any given $\kappa \in \mathcal{K}$ of sufficient regularity
with respect to $\psi$ (Theorem 2.3); (ii) necessary and
sufficient conditions exist on the function $\kappa \in
\mathcal{K}$ such that $f_{L',\kappa}, f_{L,\kappa} \in
M_+(\psi)_{\mathrm{sym},\iny}^*$ for all $L' \in BL(\NN)$ and $L
\in BL(\RR_+)$ (Theorem 2.7); and (iii) the Riesz semi-norm
$\rho_1(x)$ of $x \in M(\psi)$ is the supremum of the values
$\inset{ f_{L,\kappa}(|x|)}{L \in BL(\RR_+)}$ given certain
conditions on $\kappa$ and $\psi$ (Theorem 2.8).

\subsection{Definitions and Results}

\begin{dfnS}
Let $\psi \in \Omega_\iny$ and $\kappa \in\mathcal{K}$. Then
$\kappa$ is said to have \textbf{\emph{restricted growth with
respect to}} $\psi$ if
$$\text{F-}\lim_{n \to \iny} \frac{\psi(\kappa(n))}{\psi(\kappa(n+1))} = 1.$$
\end{dfnS}

\begin{dfnS}
Let $\psi \in \Omega_\iny$.  Denote by $\mathrm{R}(\psi)$
\textbf{\emph{the set of all $\kappa \in \mathcal{K}$ that have
restricted growth with respect to}} $\psi$.
\end{dfnS}

It is immediate the set $\mathrm{R}(\psi)$ is non-empty. The
concave function $\psi$ is an invertible function
such that $\psi^{-1}$ belongs to $\mathrm{R}(\psi) \subset \mathcal{K}$.
The rationale for introducing the set
R$(\psi)$ is provided by the following result.

\begin{thmS}
Let $\psi \in \Omega_\iny$ and $\kappa \in \mathrm{R}(\psi)$.
\begin{prop2list}{6}{2}{2}
\item Let $L \in BL(\RR_+)$.
Then there exists
$L' \in BL(\NN)$ such that
\display{f_{L,\kappa}(x) = f_{L',\kappa}(x) \ \fa x \in M_+(\psi).}
\item Let $L' \in BL(\NN)$.    Then
there exists $L \in BL(\RR_+)$ such that
\display{f_{L,\kappa}(x) = f_{L',\kappa}(x) \ \fa x \in M_+(\psi).}
\end{prop2list}
\end{thmS}

The proof of Theorem 2.3 appears in Section 2.2.
Theorem 2.3 says the sets $\inset{ f_{L,\kappa} }{ L \in
BL(\RR_+)}$ and $\inset{ f_{L',\kappa} }{ L' \in BL(\NN) }$ are
identical as sets of functionals on $M_+(\psi)$ when $\kappa \in \mathrm{R}(\psi)$.
This has an
important corollary.

\begin{corS}
Let $\psi \in \Omega_\iny$, $\kappa \in \mathrm{R}(\psi)$
and $x \in M_+(\psi)$.
Then
\display{ \text{F-}\lim_{t \to \iny} \phi_\kappa(x)(t) = A}
if and only if
\display{ \text{F-}\lim_{n \to \iny} \phi_\kappa(x)(n) = A}
for some $A \geq 0$.
\par
\preskip \proof{Immediate from Theorem 2.3.}
\end{corS}

The condition that $\kappa$ has restricted growth with respect to $\psi$
identifies the two sets of functionals as above.  However, the condition is not sufficient
to ensure \emph{additivity} of the functionals.

\begin{dfnS} We say that $\kappa\in \mathcal{K}$ is of \textbf{\emph{exponential increase}} if $\exists \, C > 0$ such that
$\forall t > 0$
\display{ \kappa(t + C) > 2 \kappa(t).}
\end{dfnS}

\begin{dfnS}
Let $\psi \in \Omega_\iny$. We denote the \textbf{\emph{set of
elements of $\mathrm{R}(\psi)$ that are of exponential increase
by}} $\mathrm{R_{exp}}(\psi)$.
\end{dfnS}

The rationale for introducing the functions of exponential
increase is provided by the following result.

\begin{thmS}Let $\psi \in \Omega_\iny$ and
$\kappa \in \mathrm{R}(\psi)$.  Then the following statements
are equivalent
\begin{prop2list}{10}{2}{4}
\item $f_{L',\kappa} \in M_+(\psi)_{\mathrm{sym},\iny}^*$
$\fa L' \in BL(\NN)$,
\item $f_{L,\kappa} \in M_+(\psi)_{\mathrm{sym},\iny}^*$
$\fa L \in BL(\RR_+)$,
\item $\kappa \in \mathrm{R_{exp}}(\psi)$.
\end{prop2list}
\par
\preskip \proof{(i) $\Leftrightarrow$ (iii) \ Let $p_n :=
\kappa(n)$ define the sequence $\{ p_n \}_{n=1}^\iny$.  Then (i)
is equivalent to: (a) the existence of $m \in \NN$ such that $2p_n
\leq p_{n+m}$, and (b) $\text{F-}\lim_{n \to \iny}
\frac{\psi(p_n)}{\psi(p_{n+1})} = 1$ by \cite[Theorem 3.8]{DPSSS}.
Since $\kappa$ is increasing, it is evident that the condition (a)
is equivalent to the assertion  that $\kappa$ is of exponential
increase.

Condition (b) is exactly the condition $\kappa$ has restricted growth
with respect to $\psi$.

(i) $\Leftrightarrow$ (ii) is immediate from Theorem 2.3. }
\end{thmS}

Theorem 2.3 and Theorem 2.7 will allow us, in following sections,
to apply the results on singular symmetric functionals in
\cite{DPSS} and \cite{DPSSS}
to the construction of Connes in \cite{CN}.
One of the results that we shall apply,
the following and final result for this section,
is a more precise version of \cite{DPSSS} Theorem 4.1.

\begin{thmS}
Let $\psi \in \Omega_\iny$. Let $\kappa\in \mathcal{K}$ be
such that
$$\lim_{n \to \iny} \frac{\psi(\kappa(n))}{\psi(\kappa(n+1))} = 1.$$
Then
$$\rho_1(x) = \sup \inset{f_{L,\kappa}(|x|)}{L \in BL(\RR_+)} \quad \forall x\in M(\psi).$$
\par
\preskip
\proof{
Let $\kappa : [0,\iny) \to [0,\iny)$ be continuous, unbounded,
and increasing such that
$$
\lim_{n \to \iny} \frac{\psi(\kappa(n))}{\psi(\kappa(n+1))} = 1
\eqno{(2.1)}
$$
Let $x \in M(\psi)$.  Without loss of generality, we assume
$\rho_1(x) =1$. Clearly
$$
q(x) := \sup \inset{L(\phi_{\kappa}(|x|))}{L \in BL(C_b)} \leq
    \limsup_{t \to \iny} \phi_{\kappa}(|x|)(t) = \rho_1(x) = 1 \eqno{(2.2)}
$$
Let
$$
a_{n} := \phi_{\kappa}(x)(n) = \frac{1}{\psi(\kappa(n))}
\int_0^{\kappa(n)} x^*(s)ds.
$$
By (2.2) there exists increasing sequence $t_k \to \iny$ such that
$$
\lim_{k \to \iny} \phi_{\kappa}(|x|)(t_k) = 1 \eqno{(2.3)}
$$
Let $n_k \in \NN$ such that $\kappa(n_k-1) \leq t_k \leq
\kappa(n_k)$. Then,
$$
\phi_{\kappa}(|x|)(t_k) \leq \frac{\psi(\kappa(n_k))}{\psi(\kappa(n_k-1))}
\frac{1}{\psi(\kappa(n_k))} \int_0^{\kappa(n_k)} x^*(s)ds
= a_{n_k} \frac{\psi(\kappa(n_k))}{\psi(\kappa(n_k-1))}.
$$
Hence
$$
a_{n_k} \geq \phi_{\kappa}(|x|)(t_k) \frac{\psi(\kappa(n_k-1))}{\psi(\kappa(n_k))}.
$$
Let $\epsilon > 0$, then by (2.1) and (2.3) there exists $K$ such
that $\fa k > K$,
$$
a_{n_k} > 1 - \frac{\epsilon}{3}.
$$
Now, for $i=1,2,...$ let $k_i$ be the smallest integer greater
than $k_{i-1}$ and $K$ such that
$$
\frac{\psi(\kappa(n_{k_i}+i))}{\psi(\kappa(n_{k_i}))} < 1 + \frac{\epsilon}{3}.
$$
The integer $k_i$ exists for each $i$ by equation (2.1). Hence for
all $j = 1,...,i$
$$
a_{n_{k_i}+j} \geq \frac{\psi(\kappa(n_{k_i}))}{\psi(\kappa(n_{k_i}+i))}
\frac{1}{\psi(\kappa(n_{k_i}))} \int_0^{\kappa(n_{k_i})} x^*(s)ds
> \frac{1}{1+\epsilon/3} a_{n_{k_i}} > \frac{1 - \epsilon/3}{1+\epsilon/3}
> 1 - \epsilon.
$$
Then, applying Sucheston's Theorem \cite{S}
and using Theorem 2.3 above, we obtain the existence of $L \in
BL(\RR_+)$ such that $L(\phi_\kappa(x)) \geq 1 - \epsilon$.  Hence
$q(x) \geq 1 - \epsilon$ for arbitrary $\epsilon > 0$ and $q(x) =
\rho_1(x) = 1$. }
\end{thmS}

\subsection{Technical Results}

This section culminates in the proof of Theorem 2.3.

\begin{lemmaS}
Let $L \in BL(\RR_+)$.  Then
\display{L'(\alpha) := L(\bp(\alpha)) \ \fa \alpha \in \ell^\iny(\NN)}
defines an element of $BL(\NN)$.  \par
\preskip
\proof{Let $L \in BL(\RR_+)$.
Then $L'$ is linear as $L$ and $p$ are linear,
$\nm{L'(\alpha)} \leq \nm{L}\nm{\bp(\alpha)} = \nm{\alpha}$
and $L'(1) = 1$ by Lemma 1.2.
Let $k \in \NN$.  Then
$L'(T_k(\alpha)) = L(\bp(T_k(\alpha))) = L(T_k(\bp(\alpha))) = L(\bp(\alpha)) = L'(\alpha)$
by Lemma 1.2(iii) and translation invariance of $L$.}
\end{lemmaS}

\begin{lemmaS}
Let $L' \in BL(\NN)$.    Then
\display{L(g) := L'(E_{\NN}(g)) \ \  \fa g \in C_b([0,\iny))}
defines an element of $BL(\RR_+)$.
\par
\preskip
\proof{
Let $L' \in BL(\NN)$.
Then $L$ is linear as $L'$ and $E_{\NN}$ is linear,
$\nm{L(g)} \leq \nm{L'}\nm{E_{\NN}(g)} \leq \nm{g}$
and $L(1) = L'(E_{\NN}(1))  = L'(1) = 1$ by Lemma 1.3.
It remains to be shown $L$ is translation invariant.
Let $a \in (0,1)$. Then
$$
\begin{array}{rl}
L(T_a(g)) = & L'(E_{\NN}(T_a(g))) = L'(\{\int_{n-1}^ng(s+a)ds\}_{n=1}^\infty)
= L'(\{\int_{n-1+a}^{n+a}g(s)ds\}_{n=1}^\infty) \\
= & L'(\{\int_{n-1+a}^{n}g(s)ds\}_{n=1}^\infty)+L'(\{\int_{n}^{n+a}g(s)ds\}_{n=1}^\infty)\\
= & L'(\{\int_{n+a}^{n+1}g(s)ds\}_{n=1}^\infty)+L'(\{\int_{n}^{n+a}g(s)ds\}_{n=1}^\infty)\\
= & L'(\{\int_n^{n+1}g(s)ds\}_{n=1}^\infty) = L'(T_1(E_{\NN}(g)))=L(g).
\end{array}
$$
Let $[b]$ be the greatest integer less than $b > 0$.
Then $T_b=T_{[b]}+T_a$ where $0\le a=b-[b]<1$.
The translation invariance of $L$ in the general case
follows from Lemma 1.3(iii).
}
\end{lemmaS}

\begin{lemmaS}
Let $\psi \in \Omega_\iny$ and $\kappa \in \mathrm{R}(\psi)$. Let
$x \in M_+(\psi)$ and define
\display{j_n(x) = \frac{1}{\psi(\kappa(n))} \int_{\kappa(n)}^{\kappa(n+1)}x^*(s)ds}
and
\display{K_n(x) = \sup_{t \in [n,n+1]} \Big| \phi_\kappa(x)(t)
    - \phi_\kappa(x)(n) \Big|.}
Then
\display{
L'( \{ j_n(x) \}_{n=1}^\iny ) = 0 = L'( \{ K_n(x) \}_{n=1}^\iny ) \ \fa L' \in BL(\NN).}
\par
\preskip
\proof{Let $x \in M_+(\psi)$.
We abbreviate notation by setting $g(t) := \phi_\kappa(x)(t)$
and $\alpha_n = \frac{\psi(\kappa(n))}{\psi(\kappa(n+1))}$.
Let $h_n(x) = g(n+1) - \alpha_n g(n)$.  Note $0 \leq \alpha_n \leq 1 \fa n$
as $\psi \circ \kappa$ is increasing.
Then $\text{F-}\lim_n \alpha_n = 1$ and $\text{F-}\lim_n |1-\alpha_n| = 0$
by hypothesis on $\kappa$.
Hence $L'(h_n(x)) = L'(T_1 \bd{g}) - (\text{F-}\lim_n \alpha_n )L'(\bd{g}) =
L'(\bd{g}) - L'(\bd{g}) = 0$
by \cite{DPSSS} Lemma 3.4 and translation invariance of $L'$.  Moreover
\display{L'(j_n(x)) = 1 \cdot
L'(j_n(x)) = (\text{F-}\lim_n \alpha_n) L'(j_n(x))
= L'(\alpha_n j_n(x)) = L'(h_n(x)) = 0}
again by \cite{DPSSS} Lemma 3.4.  Now
\begin{eqnarray*}
 K_n(x) & = & \sup_{t \in [n,n+1]} \Big|\frac1{\psi(\kappa(t))}
\int_0^{\kappa(n)}x^*ds+\frac1{\psi(\kappa(t))}
\int_{\kappa(n)}^{\kappa(t)}x^*ds-
\frac1{\psi(\kappa(n))}\int_0^{\kappa(n)}x^*ds
     \Big| \\
& \le & \sup_{t \in [n,n+1]} \Big|\phi_\kappa(x)(n) \Big( \frac{\psi(\kappa(n))}
{\psi(\kappa(t))}-1 \Big) \Big|+ \sup_{t \in [n,n+1]}
\Big|
\frac1{\psi(\kappa(t))}
\int_{\kappa(n)}^{\kappa(t)}x^*ds
\Big|  \\
& \le & \Big|\phi_\kappa(x)(n) \Big( \frac{\psi(\kappa(n))}
{\psi(\kappa(n+1))}-1 \Big) \Big|+
\Big|
\frac1{\psi(\kappa(n))}
\int_{\kappa(n)}^{\kappa(n+1)}x^*ds
\Big| \\
& \le & \nm{x}_{M(\psi)} | 1-\alpha_n | + j_n(x)
\end{eqnarray*}
Hence $L'( \{ K_n(x) \}) \leq L'( \{j_n(x) \}) = 0$ by results above.
}
\end{lemmaS}

\begin{propS}
\ Let $\psi \in \Omega_\iny$ and $\kappa \in \mathrm{R}(\psi)$.
Then $\phi_\kappa(x)$ is almost piecewise linear at infinity
for all $x \in M_+(\psi)$.
\par
\preskip \proof{We abbreviate the notation by setting $g :=
\phi_\kappa(x)$, $C_n = \sup_{t \in [n,n+1]} |g(t) - \bp
\bd{g}(t)|$ and $K_n = \sup_{t \in [n,n+1]} |g(t) - g(n)|$,
$n=0,1,2,...$ Let $f = \bp(\{ C_{n-1} \}_{n=1}^\iny)$. Then $|g -
\bp \bd{g}| \leq 2f(t+1/2)$.  Hence $L(|g - \bp \bd{g}|) \leq 2
L(T_{1/2}f) = 2 L(f)$. We now evaluate $L(f)$.  By Lemma 2.9 $L(f)
= L'(\{ C_n \}_{n=1}^\iny)$ for some $L' \in BL(\NN)$.  Consider
\vspace*{-0.5cm} \display{\begin{array}{rl} C_n = & \sup_{t \in
[n,n+1)} \Big| g(t) - \Big( g(n)
    + (g(n+1) - g(n))(t-n) \Big) \Big| \\
\leq & \sup_{t \in [n,n+1)} | g(t) - g(n)|
    + |g(n+1) - g(n)| \ \leq \ 2 K_n.
\end{array}}
Hence $L(f) = L'(\{ C_n \}_{n=1}^\iny) \leq 2
L'(\{K_n\}_{n=1}^\iny) = 0$ by Lemma 2.11.}
\end{propS}

We are now in a position to prove Theorem 2.3.

\subsubsection*{Proof of Theorem 2.3}

\noindent (i) \ Let $L \in BL(\RR_+)$ and $L' \in BL(\NN)$ as in
Lemma 2.9. By Proposition 2.12 $L(\phi_\kappa(x) - \bp
\bd{\phi_\kappa(x)}) = 0$, hence $L(\phi_\kappa(x)) = L(\bp
\bd{\phi_\kappa(x)}) = L'(\bd{\phi_\kappa(x)})$.

\medskip \noindent (ii) \ Let $L' \in BL(\NN)$.
Let $L \in BL(\RR_+)$ as in Lemma 2.10, $L(\phi_\kappa(x)) =
L'(E_{\NN}(\phi_\kappa(x)))$. Since
$E_{\NN}(\phi_\kappa(x))(n)=\phi_\kappa(x)(\xi_n)$ for some
$\xi_n\in[n-1,n]$ for each $n \in \NN$ then \display{
|\phi_\kappa(x)(n)-E_{\NN}(\phi_\kappa(x))(n)|\le\sup_{t\in[n-1,n]}|\phi_\kappa(x)(t)-\phi_\kappa(x)(n)|=K_n(x)}
for each $n \in \NN$.  Consequently, by Lemma 2.11,
$$
\begin{array}{rcl}
|L'(\bd{\phi_\kappa(x)}) - L'(E_{\NN}(\phi_\kappa(x)))| & \le
&  L'(\{|\phi_\kappa(x)(n)-E_{\NN}(\phi_\kappa(x))(n)|\}_{n=1}^\infty) \\
& \le & L'(K_n(x)) = 0
\end{array}
$$
or $L(\phi_\kappa(x))=L'(\bd{\phi_\kappa(x)})$ as required. \quad
$\Box$

\section{Measurability in Marcinkiewicz Spaces}

\noindent Let $\psi \in \Omega_\iny$ and $\kappa \in \mathcal{K}$.
Having constructed the family $\inset{f_{L,\kappa}}{L \in
BL(\RR_+)}$ of functionals on $M_+(\psi)$, it is natural to
consider elements $x$ in $M_+(\psi)$ such that
$$
f_{L_1,\kappa}(x) = f_{L_2,\kappa}(x) \fa L_1,L_2 \in BL(\RR_+).
$$
It is obvious that the equation above holds if and only if
$\phi(x)(t)$ is almost convergent (see Definition 1.1) and
consequently, in this case
$$
f_{L,\kappa}(x) = \text{F-} \lim_{t \to \iny} \phi_\kappa(x)(t) = A
$$
for some $A \geq 0$ and all $L \in BL(\RR_+)$. Necessary and
sufficient conditions for almost convergence, even for sequences
\cite{GL}, are somewhat complicated . In studying the function
$\phi_\kappa(x)$ it is more preferable to consider the notions of
Cesaro convergence (definition below) and ordinary convergence and
`squeeze' almost convergence inbetween. In this section we: (i)
establish Cesaro convergence is weaker than almost convergence
which in turn is weaker than ordinary convergence (Remark 3.1,
Corollary 3.4), and then (ii) consider Tauberian conditions (see
\cite[Section 6.1]{Hardy}) on the function $\phi_\kappa(x)$ such
that Cesaro convergence implies ordinary convergence and hence the
notions of Cesaro, almost and ordinary convergence are identical
for $\phi_\kappa(x)$ (Theorem 3.7, Corollary 3.10).

\subsection{Definitions and Results}

Let $\{a_n \}_{n=1}^\iny \in \ell^\iny(\NN)$.  Define
\display{b_n(p) = \frac{1}{n} \sum_{i=0}^{n-1} a_{p+i}}
for $p \in \NN$.   We recall
from \cite{GL} that $\{ a_n \}_{n=1}^\iny$ is almost convergent
(see Definition 1.1) if and only if
\display{L'(\{ a_n \}) = \text{F-}\lim_n a_n = \lim_n b_n(p) = A}
for all $L' \in BL(\NN)$ where $\lim_n b_n(p) = A$ uniformly with respect to $p \in \NN$.

\medskip \noindent A \textbf{sequence} $\{ a_n \}_{n=1}^\iny$ is called
\begin{prop2list}{16}{2}{4}
\item \textbf{Cesaro convergent} if $\lim_n b_n(1) = A$
\item \textbf{almost convergent} if $\lim_n b_n(p) = A$
uniformly with respect to $p \in \NN$
\item \textbf{convergent} if $\lim_n a_n = A$
\end{prop2list}
for some $A \geq 0$. We denote by C, F and S \textbf{the sets of
all Cesaro convergent sequences, almost convergent sequences and
convergent sequences}, respectively.

\REMS{Since \display{\lim_n a_n = A \ \implies \ \text{F-}\lim_n
a_n = \lim_n b_n(p) = A \ \implies \ \lim_n b_n(1) = A} we have
the inclusion of sets ${ \rm S \subset \rm F \subset \rm
C.}$ }

\vspace*{0cm}

\begin{dfnS}
Let $\psi \in \Omega_\iny$ and $\kappa \in\mathcal{K}$. Let $x \in
M_+(\psi)$.  We say $x$ is
\begin{prop2list}{16}{2}{4}
\item \emph{C}$_\kappa$\textbf{\emph{-measurable}} if
$\bd{\phi_\kappa(x)} \in \rm C$, \item \emph{F}$_\kappa$
\textbf{\emph{-measurable}} if $\bd{\phi_\kappa(x)} \in \rm F$,
\item \emph{S}$_\kappa$\textbf{\emph{-measurable}} if
$\bd{\phi_\kappa(x)} \in \rm S$.
\end{prop2list}
\end{dfnS}

\bigskip \noindent Define for $\mu > 0$, \display{C(g)(\mu) = \frac{1}{\mu}
\int_{0}^\mu g(t) dt.} A \textbf{function} $g\in C_b([0,\iny))$ is called
\begin{prop2list}{16}{2}{4}
\item \textbf{Cesaro convergent} if $\lim_{t \to \iny} C(g)(t) = A$
\item \textbf{almost convergent at infinity} if
F-$\lim_{t \to \iny} g(t) = A$
\item \textbf{convergent at infinity} if $\lim_{t \to \iny} g(t) = A$
\end{prop2list}
for some $A \geq 0$. We denote by $\mathcal{C}$, $\mathcal{F}$ and
$\mathcal{S}$ \textbf{the sets of all Cesaro convergent functions,
almost convergent functions and functions convergent at infinity},
respectively.

\vspace*{0cm}
\begin{thmS}
Let $g \in C_b([0,\iny))$.  Then
\display{ [a,b] \subset \inset{L(g)}{L \in BL(\RR_+)} }
where
\display{\quad a=\liminf_{t\to\iny}C(g)(t),\
b=\limsup_{t\to\iny}C(g)(t).}
\par
\preskip \proof{Suppose the result is false.  Then there exists $c
\in [a,b]$ such that $c \not= L(g)$ for any $L \in BL(\RR_+)$.  By
continuity of $C(g)$ there exists a sequence $t_n\to\iny$ as $n
\to \iny$ such that $C(g)(t_n) \to c$. Let us consider
$C_b([0,\iny))^*$ equipped with the weak$^*$-topology. Then the
unit ball $B$ of $C_b([0,\iny))^*$ is weak$^*$-compact. Hence, the
sequence of functionals $\delta_{t_n}(f)=f(t_n),\ n=1,2,...,$ has
a limit point $V \in B$. In fact, this limit point belongs to the
weak$^*$ compact subset $B_1$ of positive elements $\gamma$ of the
unit ball  $B\subset C_b([0,\iny))^*$ such that $\gamma(1)=1$.
From weak$^*$ convergence the state $V$ has the following
properties, (i) $V(p)=\lim_n p(t_n)=0$ for every function $p \in
C_0([0,\iny))$, and (ii) $V(C(g))=\lim_n C(g)(t_n)=c$.

Define the functional $L(f):=V(C(f))$ for $f \in C_b([0,\iny))$.
It is immediate that $L(g) = c$ by property (ii).  Hence,
if $L$ belongs to $BL(\RR_+)$, the supposition on $c$ is false
and the result is proven.

We show the functional $L$ is translation
invariant. Indeed, for any $f\in C_b([0,\iny))$
$$
C(T_af)(\mu)-C(f)(\mu)=\frac1\mu\int_0^\mu [f(t+a)-f(t)]dt=
\frac1\mu\left[\int_0^a f(t)dt +\int_\mu^{\mu+a}f(t)dt\right]
\to 0
$$
for $\mu\to\infty$. Hence
translation invariance of $L$ follows by property (i).
Trivially $L(1) = V(C(1)) = V(1) = 1$.  Hence
$L \in BL(\RR_+)$.
}
\end{thmS}

\begin{corS}
Let $\mathcal{C}$, $\mathcal{F}$ and $\mathcal{S}$ be the sets defined
as above.  Then
\display{\mathcal{S} \subset \mathcal{F} \subset \mathcal{C}.}
\proof{The inclusion $\mathcal{S} \subset \mathcal{F}$ is immediate.
The inclusion $\mathcal{F} \subset \mathcal{C}$ is immediate
from Theorem 3.3 }
\end{corS}

\begin{dfnS}
Let $\psi \in \Omega_\iny$ and $\kappa \in \mathcal{K}$. Let $x
\in M_+(\psi)$.  We say $x$ is
\begin{prop2list}{16}{2}{4}
\item $\mathcal{C}_\kappa$\textbf{\emph{-measurable}} if
$\phi_\kappa(x) \in \mathcal{C}$, \item
$\mathcal{F}_\kappa$\textbf{\emph{-measurable}} if $\phi_\kappa(x)
\in \mathcal{F}$, \item
$\mathcal{S}_\kappa$\textbf{\emph{-measurable}} if $\phi_\kappa(x)
\in \mathcal{S}$, \item $\mathcal{S}$\textbf{\emph{-measurable}}
if $\phi(x) \in \mathcal{S}$.
\end{prop2list}
\end{dfnS}
\vspace*{0cm} \REMS{We draw the  reader's attention to the fact
that since $\kappa$ is continuous, $x$ is
$\mathcal{S}_\kappa$-measurable  if and only if  $x$ is
$\mathcal{S}$-measurable.  The same (simple) analysis does not
work with the notion of $\rm S_\kappa$-measurability introduced in
Definition 3.2. Nevertheless, it is established in the following
theorem that the equivalence of ${\rm S}_\kappa$-measurability of
an element $x$ with $\mathcal{S}$-measurability of $x$ holds under
natural restriction on $x$.}

\begin{thmS}
Let $\psi \in \Omega_\iny$ and $\kappa \in \mathrm{R}(\psi)$. Let
$x \in M_+(\psi)$ be such that
$$
t \, \phi_\kappa(x)'(t) > -H
$$
for some $H>0$ and all $t > 0$.  Then the following statements are
equivalent
\begin{prop2list}{16}{2}{3}
\item $x$ is $\rm C_\kappa$-measurable, \item $x$ is
$\mathcal{C}_\kappa$-measurable, \item $x$ is $\rm
F_\kappa$-measurable, \item $x$ is
$\mathcal{F}_\kappa$-measurable, \item $x$ is $\rm
S_\kappa$-measurable, \item $x$ is $\mathcal{S}$-measurable.
\end{prop2list}
\end{thmS}

The proof of Theorem 3.7 appears in Section 3.2.
The hypothesis on the derivative $\phi_\kappa(x)'$, which depends on
$x\in M_+(\psi)$,
can be made independent of $x$ by a stronger hypothesis on the function $\kappa$.
We recall that $\kappa \in \mathcal{K}$ is an invertible differentiable
function such that $\kappa(0) = 0$.

\begin{dfnS}
Let $\psi \in \Omega_\iny$ and $\kappa \in \mathcal{K}$. We say
$\kappa$ has \textbf{\emph{dominated growth with respect to}}
$\psi$ if $\exists C > 0$ such that $\forall t > 0$
$$
\frac{(\psi\circ\kappa)'(t)}{\psi\circ\kappa(t)} < \frac{C}{t}.
$$
Denote by $\mathrm{D}(\psi)$ \textbf{\emph{the set of $\kappa \in
\mathcal{K}$ that have dominated growth with respect to $\psi$}}.
\end{dfnS}

\noindent It is immediate that the set $\mathrm{D}(\psi)$ is non-empty
since it contains $\psi^{-1}$.
The rationale for introducing the set
D$(\psi)$ is provided by the following result.

\begin{corS}
Let $\psi \in \Omega_\iny$, $\kappa \in \mathrm{D}(\psi)$ and $x
\in M_+(\psi)$. Then the following statements are equivalent
\begin{prop2list}{16}{2}{3}
\item $x$ is $\rm C_\kappa$-measurable, \item $x$ is
$\mathcal{C}_\kappa$-measurable, \item $x$ is $\rm
F_\kappa$-measurable, \item $x$ is
$\mathcal{F}_\kappa$-measurable, \item $x$ is $\rm
S_\kappa$-measurable, \item $x$ is $\mathcal{S}$-measurable.
\end{prop2list}
\end{corS}

The proof of Corollary 3.9 also appears in Section 3.2.
In terms of the functionals $f_{L,\kappa}$ of Definition 1.7
the preceding result may be reformulated as follows.

\begin{corS}
Let $\psi \in \Omega_\iny$, $\kappa \in \mathrm{D}(\psi)$
and $x \in M_+(\psi)$.
Then the following statements are equivalent
\begin{prop2list}{10}{2}{4}
\item $x$ is $\mathcal{C}_\kappa$-measurable,
\item $f_{L,\kappa}(x)$ is independent of $L \in BL(\RR_+)$,
\item $f_{L,\kappa}(x) = \lim_{t \to \iny} \phi(x)(t) \ \fa L \in BL(\RR_+).$
\end{prop2list}
\end{corS}

The equivalence of statements (ii) and (iii) in the above
Corollary is a new and surprising result.  The implication
of the result may be seen
in the context of the work of A. Connes.
For this end we introduce notions relevant to \cite{CN}.

\begin{dfnS}
Let $\psi \in \Omega_\iny$ and $\kappa \in\mathcal{K}$.  Then we
say $x \in M_+(\psi)$ is
\begin{prop2list}{10}{2}{4}
\item $\kappa$\textbf{\emph{-measurable}} if $f_{L,\kappa}(x)$ is
independent of $L \in BL(\RR_+)$, and \item
\textbf{\emph{Tauberian}} if $$\lim_{t \to \iny}\phi(x)(t) =
\lim_{t \to \iny} \frac{1}{\psi(t)} \int_0^t x^*(s)ds = A$$ for
some $A \geq 0$.
\end{prop2list}
\end{dfnS}

\begin{dfnS}
Let $\psi \in \Omega_\iny$. Denote by ${\cal M}_\kappa^+(\psi)$
(respectively, $\mathcal{T}_+(\psi)$) \textbf{\emph{the set of
$\kappa$-measurable}} (respectively, \textbf{\emph{Tauberian}})
\textbf{\emph{elements of $M_+(\psi)$}}. We also define the set
${\cal M}_+(\psi) := \cap_{\kappa \in \mathrm{R_{exp}}(\psi)}
{\cal M}_\kappa^+(\psi)$ called \textbf{\emph{the set of
measurable positive elements of the Marcinkiewicz space}}
$M(\psi)$.
\end{dfnS}
\begin{thmS}
Let ${\cal M}_\kappa^+(\psi)$ and
${\cal M}_+(\psi)$ be defined as above.
Then
\begin{prop2list}{10}{2}{4}
\item ${\cal M}_\kappa^+(\psi)$ is a closed, symmetric
subcone of $M_+(\psi)$ when $\kappa \in \mathrm{R_{exp}}(\psi)$,
\item ${\cal M}_+(\psi)$ is a closed, symmetric
subcone of $M_+(\psi)$.
\end{prop2list}
\par
\preskip \proof{(i) Closedness, symmetricity and additivity follow
from the fact $f_{L,\kappa}$ is an additive singular symmetric
functional on $M_+(\psi)$ by Theorem 2.7.  (ii)  Follows from (i)
as ${\cal M}_+(\psi)= \cap_{\kappa \in \mathrm{R_{exp}}} {\cal
M}_\kappa^+(\psi)$. }
\end{thmS}

\medskip \noindent The implication of Theorem 3.7 is the following result
which connects
Proposition IV.2.$\beta$.4 and Proposition IV.2.$\beta$.6 of \cite{CN}.
We shall elaborate on this result in Section 5 and the implications
of the result for non-commutative geometry in the concluding section.

\begin{thmS}
Let $\psi \in \Omega_\iny$. Then
$$
\mathcal{T}_+(\psi) = \mathcal{M}^+_\kappa(\psi)
$$
for all $\kappa \in \mathrm{D}(\psi)$ and, if there exists
$\kappa \in \mathrm{D}(\psi)$ of exponential increase,
$$
\mathcal{T}_+(\psi) = \mathcal{M}_+(\psi)
             = \mathcal{M}^+_\kappa(\psi)
$$
for all $\kappa \in \mathrm{D}(\psi)$.
\par
\preskip \proof{The first result is immediate from Corollary 3.9.
Suppose $\kappa_1 \in \mathrm{D}(\psi)$ is of exponential
increase. Then $\kappa_1 \in \mathrm{R_{exp}}(\psi)$ by
Proposition 3.20(i) of next section. Hence $\mathcal{T}_+(\psi)
\subset \mathcal{M}_+(\psi) \subset
 \mathcal{M}^+_{\kappa_1}(\psi) = \mathcal{T}_+(\psi) =
\mathcal{M}^+_{\kappa}(\psi)$ for any $\kappa \in \mathrm{D}(\psi)$,
where the last equality is given by the first result.}
\end{thmS}

\REMS{ It was shown in \cite{DPSSS} that if
$$
\liminf_{t\to\iny} \frac{\psi(2t)}{\psi(t)}=1 \ \ \mbox{but} \ \
\limsup_{t\to\iny} \frac{\psi(2t)}{\psi(t)}=2
$$
then there exists $x_0 > 0 $ in $M(\psi)\setminus M_1(\psi)$ such
that all additive symmetric functionals defined on $M(\psi)$
vanish on $x_0$.  However, if $\phi(x)(t) \to 0$ as $t \to \iny$
then $\rho_1(x_0) = 0$ and $x_0 \in M_1(\psi)$, which is a
contradiction.  This example shows the set $\mathcal{M}_+(\psi)$
of measurable elements and the set $\mathcal{T}_+(\psi)$ of
Tauberian elements are not the same in general and the set
$\mathrm{D}(\psi)$ can fail to admit an element of exponential
increase.  Necessary and sufficient conditions on the concave
function $\psi$ such that $\mathrm{D}(\psi)$ admits an element of
exponential increase are given in Proposition 3.20 of next
section.}

\subsection{Technical Results}

Let $\{ a_n \}_{n \in \NN} \subset \RR$ be a sequence and
$s_n = \sum_{m=1}^n a_m$ denote
the $n^{\mathrm{th}}$-partial sum.
Hardy's section on Tauberian theorems for Cesaro summability \cite[Section 6.1]{Hardy}
contains the following result.

\medskip \noindent \textbf{THEOREM 64}  If $\lim_{n \to \iny} \frac{1}{n} \sum_{m=1}^n s_m = A$
and $na_n > -H$ for some $A \in \RR$ and $H > 0$, then $\lim_{n \to \iny} s_n = A$.

\medskip \noindent We recall that any sequence $\{ b_n \}_{n \in \NN}$ is the
sequence of partial sums of the sequence $\{ a_n := b_n - b_{n-1} \}_{n\in \NN}$ with the
convention $b_0=0$.  Hence a trivial corollary of Theorem 64 is the following.

\begin{thmS}  Let $\{ b_n \}_{n \in \NN}$ be a sequence such that $b_n \geq 0$
and $n(b_n - b_{n-1}) > -H$ for some $H > 0$.
Then $\lim_{n \to \iny} \frac{1}{n} \sum_{m=1}^n b_m = A$ for some $A \geq 0$
if and only if $\lim_{n \to \iny} b_n = A$.
\end{thmS}

\noindent A continuous analogy of Theorem 64 exists in \cite[Section 6.8]{Hardy}.  It has
the following corollary.

\begin{thmS}  Let $b(t)$ be a positive piecewise differentiable function
such that $tb'(t) > -H$ for some $H > 0$ and almost all $t > 0$.
Then $\lim_{t \to \iny} \frac{1}{t} \int_0^t b(s)ds = A$ for some
$A \geq 0$ if and only if
$\lim_{t \to \iny} b(t) = A$.
\end{thmS}

These theorems are sufficient to prove Theorem 3.7 with the following
lemma.

\begin{lemmaS}
Let $b(t)$ be a piecewise differentiable function such that
$tb'(t) > -H$ for some $H > 0$ and almost all $t > 0$. Then
$n(b(n) - b(n-1)) > -2H$ for all $n \in \NN$. \par \preskip
\proof{Let $n \in \NN$.  Then $b(n) - b(n-1) \geq \inf_{t \in
[n-1,n]} b'(t) > \inf_{t \in [n-1,n]} -Ht^{-1} \geq - H(n-1)^{-1}
\geq -2Hn^{-1}$.}
\end{lemmaS}

\subsubsection*{Proof of Theorem 3.7}  The scheme of implications
shall be \\
\centerline{ \begin{tabular}[r]{c@{}c@{}c@{}c@{}c}
\text{ (i) } & $\Leftrightarrow$ & \text{ (iii) } &
$\Leftrightarrow$ & \text{ (v) } \\[-6pt]
  & & $\Updownarrow$ & & \\[-6pt]
 \text{ (ii) } & $\Leftrightarrow$ & \text{ (iv) } &
 $\Leftrightarrow$ & \text{ (vi). }
\end{tabular} } \\

\noindent \text{(i)} $\Leftarrow$ \text{(iii)} $\Leftarrow$
\text{(v)} is Remark 3.1
and (v) $\Leftarrow$ (i) is provided by
Lemma 3.18 and Theorem 3.16 using $b(t) = \phi_\kappa(t)$.

\noindent \text{(iii)} $\Leftrightarrow$ \text{(iv)} is Corollary 2.4.

\noindent \text{(ii)} $\Leftarrow$ \text{(iv)} $\Leftarrow$
\text{(vi)} is Corollary 3.4 and (vi) $\Leftarrow$ (ii) is
provided by Theorem 3.17 using $b(t) = \phi_\kappa(t)$ and
Remark 3.6. \vspace*{-.5cm}
\begin{flushright} $\Box$ \end{flushright}

The following Propositions are sufficient to prove Corollary 3.9.

\begin{propS}
Let $\psi \in \Omega_\iny$ and $\kappa \in \mathrm{D}(\psi)$. Then
$t\phi_\kappa(x)'(t) > - C \nm{x}_{M(\psi)}$ for all $t > 0$.
\par
\preskip \proof{From the proof of Proposition 1.8
$$
\phi(x)'(t) \geq - \frac{\psi'(t)}{\psi(t)} \nm{x}_{M(\psi)}
$$
The substitution $t \mapsto \kappa(t)$, multiplication of both
sides by the positive number $t\kappa'(t)$ for $t > 0$ and the
elementary property $(f\circ\kappa)'(t) = f'(\kappa(t))\kappa'(t)$
yields
$$
t\phi_\kappa(x)'(t) \geq
 - t\frac{(\psi \circ \kappa)'(t)}{\psi \circ \kappa(t)} \nm{x}_{M(\psi)}.
$$
The result now follows from the hypothesis $\kappa \in
\mathrm{D}(\psi)$. }
\end{propS}

\begin{propS}Let $\psi \in \Omega_\infty$.  Then
\begin{prop2list}{10}{2}{2}
\item $\mathrm{D}(\psi) \subset \mathrm{R}(\psi)$
\end{prop2list}
and the following statements are equivalent
\begin{prop2list}{10}{2}{2}
\item[\text{(ii)}] the set $\mathrm{D}(\psi)$ contains an element
$\kappa$ of exponential increase \emph{;} \item[\text{(iii)}] $\exists \,
C>0$ such that
$$
\frac{\psi(2t)}{\psi(t)}=1+O\left(\frac1{\psi(t)^{1/{C}}}\right);
$$
\item[\text{(iv)}]  $\psi^{-1}(t^C)$ is of exponential increase
for some $C \geq 1$.

\end{prop2list}
\par
\preskip \proof{(i) Let $\kappa \in \mathrm{D}(\psi)$. Then by
Definition 3.8
$$
\log\left(\frac{(\psi \circ \kappa)(t+T)}{(\psi \circ
\kappa)(t)}\right)= \int_t^{t+T}\frac{(\psi \circ
\kappa)'(s)}{(\psi \circ \kappa)(s)}ds <
C\int_t^{t+T}s^{-1}ds=C\log\frac{(t+T)}{t}.
$$
Consequently
$$
\frac{(\psi \circ \kappa)(t+T)}{(\psi \circ \kappa)(t)}
<\left(\frac{t+T}{t}\right)^C=1+O(t^{-1}) \mbox{ for large }t.
\eqno{(3.1)}
$$
Taking $t=n$ and $T=1$ we get (i).

(ii) $\Rightarrow$ (iii) Substituting $t=1$ and $T=u-1$ into (3.1)
we get
$$
(\psi \circ \kappa)(u)<(\psi \circ \kappa)(1)u^C=O(u^C)
\eqno{(3.2)}
$$
Then, taking $T=D$ where $D$ is such that $\kappa(t+D)>2\kappa(t)$
for all $t>0$,
$$
1<\frac{\psi(2\kappa(t))}{\psi(\kappa(t))}< \frac{(\psi \circ
\kappa)(t+T)}{(\psi \circ \kappa)(t)}<1+O(t^{-1})<
1+O\left(\frac1{(\psi \circ \kappa)(t)^{1/C}}\right),
$$
where the last inequality follows from (3.2). We obtain the result
by the substitution $\kappa(t)\to t$.

(iii) $\Rightarrow$ (iv) Let $t=\psi^{-1}(u)$. Then for
sufficiently large $H$ we have
$$
\frac{\psi(2\psi^{-1}(u))}{u}<1+\frac{H}{(u)^{1/{C}}},
$$
or
$$
\psi(2\psi^{-1}(u))<u+Hu^\frac{C-1}{C}.
$$
Applying $\psi^{-1}(\cdot)$ to both sides of the last inequality
we get
$$
2\psi^{-1}(u)<\psi^{-1}(u+Hu^\frac{C-1}{C}).%
$$
If $C<1$ then $\psi^{-1}(u+Hu^\frac{C-1}{C})<\psi^{-1}(u+H)$ for
$u>1$ and in this case $\psi^{-1}(u)$ is of exponential increase.

\noindent If $C>1$ then replacing $u$ by $u^C$ we get
$$
2\psi^{-1}(u^C)<\psi^{-1}(u^C+Hu^{C-1})\le\psi^{-1}((u+H/C)^C).
$$
Consequently, $\psi^{-1}(u^C)$ is of exponential increase.

(iv) $\Rightarrow$ (ii) is immediate.

}
\end{propS}

\medskip We can now prove Corollary 3.9.

\subsubsection*{Proof of Corollary 3.9}
Let $\kappa \in \mathrm{D}(\psi)$.  Then $\kappa \in
\mathrm{R}(\psi)$ by Proposition 3.20(i) and for each $x \in
M_+(\psi)$ there exists $H = C\nm{x}_{M(\psi)} > 0$ such that
$t\phi_\kappa(x) > - H$ by Proposition 3.19.  Hence the conditions
of Theorem 3.7 are satisfied. \vspace*{-.5cm}
\begin{flushright} $\Box$ \end{flushright}

\section{Summary and Examples}

Let $\psi \in \Omega_\iny$.
For the convenience of the reader, we summarize the hypotheses on
$\kappa$ that have appeared in the previous sections.

\begin{dfnS}
Let $\psi \in \Omega_\iny$ and $\kappa \in\mathcal{K}$.  Then we
say $\kappa$
\begin{prop2list}{10}{2}{3}
\item has \textbf{\emph{restricted growth with respect to}} $\psi$
if \display{\text{F-}\lim_{n \to \iny}
    \frac{\psi(\kappa(n))}{\psi(\kappa(n+1))} = 1,}
and \textbf{\emph{the set of $\kappa$ with restricted growth with
respect to $\psi$ is denoted}} $\mathrm{R}(\psi)$. \item has
\textbf{\emph{strong restricted growth with respect to}} $\psi$ if
\display{\lim_{n \to \iny}
    \frac{\psi(\kappa(n))}{\psi(\kappa(n+1))} = 1,}
and \textbf{\emph{the set of $\kappa$ with strong restricted
growth with respect to $\psi$ is denoted}} $\mathrm{SR}(\psi)$.
\item has \textbf{\emph{dominated growth with respect to}} $\psi$
if $\exists \, C > 0$ such that $\forall t > 0$
\display{\frac{(\psi\circ\kappa)'(t)}{\psi\circ\kappa(t)} <
\frac{C}{t},} and \textbf{\emph{the set of $\kappa$ with dominated
growth with respect to $\psi$ is denoted}} $\mathrm{D}(\psi)$.
\item is of \textbf{\emph{exponential increase}} if $\exists \, C
> 0$ such that $\forall t > 0$ \display{ \kappa(t + C) > 2
\kappa(t),} and \textbf{\emph{the set of $\kappa$ of exponential
increase is denoted}} $\mathcal{K}_{\mathrm{exp}}$.
\end{prop2list}
\end{dfnS}

We denote $\mathrm{X_{\mathrm{exp}}}(\psi) = \mathrm{X}(\psi) \cap
\mathcal{K}_{\mathrm{exp}}(\psi)$, where X is D, SR, or R. The
conditions (i), (ii) and (iii) are increasingly stronger
conditions by Proposition 3.20, hence $\mathrm{D}(\psi) \subset
\mathrm{SR}(\psi) \subset \mathrm{R}(\psi)$.   We recall that
$\kappa \in \mathrm{R}(\psi)$ was sufficient for Theorem 2.3,
$\kappa \in \mathrm{R_{exp}}(\psi)$ was necessary and sufficient
for Theorem 2.7, $\kappa \in \mathrm{SR}(\psi)$ was sufficient for
Theorem 2.8, and $\kappa \in \mathrm{D}(\psi)$ was sufficient for
Corollary 3.9.  Hence $\kappa \in \mathrm{D_{exp}}(\psi)$ is the
strongest hypothesis on $\kappa$ and implies Theorems 2.3, 2.7,
2.8, 3.14 and Corollary 3.9.

We now point out some explicit examples of functions $\psi$ and
$\kappa$ for which $\kappa \in \mathrm{D_{exp}}(\psi)$.
Indeed, the functions given in Example
4.3 below appear in \cite{CN}. Consequently, Theorems 2.3, 2.7, 2.8,
3.14 and Corollary 3.9 apply to the functionals on Marcinkiewicz
operator spaces used in \cite{CN}. We shall elaborate on this in
our concluding section.

\EXS{Let $\psi : [0,\iny) \to [0,\iny)$ be a continuous, concave
and invertible function such that the inverse $\psi^{-1} :
[0,\iny) \to [0,\iny)$ is of exponential increase. Then $\psi^{-1}
\in \mathrm{D_{exp}}(\psi)$.}

\EXS{Define the function $\psi : [0,\iny) \to [0,\iny)$ by
\display{\psi(t) := \log(1+t).} Then $\psi$ is continuous, concave
and invertible. The function \display{\psi^{-1}(t) = e^t - 1} is
of exponential increase.  Hence Example 4.2 applies and $\psi^{-1}
\in \mathrm{D_{exp}}(\psi)$. The function ${\kappa}$
given by
\display{{\kappa}(t) := e^t}
is an element of the equivalence class $[\psi^{-1}]$ by Remark 1.9 and
hence provides the same set of functionals as $\psi^{-1}$.}

Let $\psi \in \Omega_\iny$.  We conclude the summary with a result
on the existence of the sets $\mathrm{X_{\mathrm{exp}}}(\psi)$
where X is D, SR, or R.

\begin{thmS}
Let $\psi \in \Omega_\iny$.  The following set A of statements are
equivalent \\
\begin{prop2list}{12}{2}{3}
\item[\text{A(i)}]  $\liminf_{t \to \iny} \frac{\psi(2t)}{\psi(t)}
= 1 \ ;$ \item[\text{A(ii)}]
$\mathrm{R_{exp}}(\psi)$ is non-empty. \\
\end{prop2list}
The following set B of statements are
equivalent \\
\begin{prop2list}{12}{2}{3}
\item[\text{B(i)}] $\lim_{t \to \iny} \frac{\psi(2t)}{\psi(t)} = 1
\ ;$ \item[\text{B(ii)}]
$\mathrm{SR_{exp}}(\psi)$ is non-empty. \\
\end{prop2list}
The following set C of statements are
equivalent \\
\begin{prop2list}{12}{2}{3}
\item[\text{C(i)}] $\frac{\psi(2t)}{\psi(t)} - 1 =
O(\psi(t)^{-1/C})$ for some $C > 0$ \emph{;} \item[\text{C(ii)}]
$\mathrm{D_{exp}}(\psi)$ is non-empty.
\end{prop2list}
\par
\preskip \proof{ Set A. \ Follows from \cite{DPSS} Theorem 3.4(i)
and \cite{DPSSS} Lemma 3.9 with Theorem 2.7.

\medskip \noindent Set B. \ (i) $\Rightarrow$ (ii) Let $\beta(t) := 2^t$.
It is immediate $\beta \in \mathrm{R_{exp}}(\psi)$ and $\lim_{t
\to \iny} \psi(2^{t+1}) / \psi(2^t) =1$ by hypothesis on $\psi$.

\smallskip \noindent (ii) $\Rightarrow$ (i) The hypothesis implies
for any $m \in \NN$ \display{\lim_{n \to \iny}
\frac{\psi(\kappa(n))}{\psi(\kappa(n+m))} = \lim_{n \to \iny}
\frac{\psi(\kappa(n+m-1))}{\psi(\kappa((n+m-1) +1))}...
\frac{\psi(\kappa(n))}{\psi(\kappa(n+1))} = 1.} Let $m'$ be any
integer greater than the $C > 0$ such that $\kappa(t + C) > 2
\kappa(t)$ for $t > 0$.  Then
\display{\frac{\psi(\kappa(n))}{\psi(\kappa(n+1 + m'))} \leq
\frac{\psi(\kappa(n))}{\psi(2\kappa(n + 1))} \leq
\frac{\psi(t)}{\psi(2t)} \leq 1 } for all $\kappa(n) \leq t \leq
\kappa(n+1)$.  Since $\kappa$ is of exponential increase then
$\kappa(n) \to \iny$ as $n \to \iny$.  Hence \display{1 = \lim_{n
\to \iny} \frac{\psi(\kappa(n))}{\psi(\kappa(n+1+m'))} \leq
\lim_{t \to \iny}  \frac{\psi(t)}{\psi(2t)} \leq 1.}

\medskip \noindent Set C. \ Follows from Proposition 3.20. }
\end{thmS}

\REMS{The example $\psi(t)=(\log(1+t))^C,\ C>1$ for large $t$ and
linear for small $t$ shows that the constant $1/C$ in Theorem 4.4
C(ii) cannot \mbox{be replaced with 1.}}

\section{Generalization of the Connes-Dixmier construction}

\subsection{Connes-Dixmier Functionals on Marcinkiewicz Spaces}

The Connes-Dixmier construction of \cite[IV.2]{CN}, which we shall
continue to clothe in the language of singular symmetric
functionals on Marcinkiewicz spaces until the concluding section,
generates singular symmetric functionals on $M_+(\psi)$ supported
at infinity for the specific function $\psi(t)=\log(1+t)$. We
recall the idea of A. Connes' method.

\begin{dfnS}
Let $SC_b^*([0,\infty))$ denote the set of all positive linear
functionals $\gamma$ on $C_b([0,\iny))$ such that $\gamma(1)=1$
and $\gamma(f) = 0$ for all $f$ in $C_0([0,\iny))$.
\end{dfnS}

A. Connes defines a symmetric functional supported at infinity on
the cone of positive elements of the Marcinkiewicz space
$M(\log(1+t))$ by the formula
$$\tau_\gamma(x) :=
 \gamma \left(
\frac{1}{\log(1 + \lambda)}\int_0^\lambda \phi(x)(u) d\log(1+u)
   \right)
$$
for all $x \in M_+(\log(1+t))$, where
$\gamma\in SC_b^*([0,\infty))$ and
$$
\phi(x)(u)=\frac{1}{\log(1+u)} \int_0^u x^*(s)ds.
$$

We generalize the construction to any Marcinkiewicz space
$M(\psi)$ of Lebesgue measurable functions, $\psi \in \Omega_\iny$, and
demonstrate the functionals so constructed are a sub-class of
functionals of the form $f_{L,\kappa}$ already studied in this
paper.

Let $k \in \mathcal{K}$.
Define
$$
M_k(g)(\lambda) := \frac{1}{k(\lambda)}\int_0^\lambda g(s) dk(s)
$$
where $g \in C_b([0,\iny))$ and $\lambda > 0$.

\begin{dfnS}
Let $\psi \in \Omega_\iny$ and $k \in\mathcal{K}$. Let $\gamma \in
SC_b^*([0,\iny))$. Then \display{\tau_{\gamma,k}(x) = \gamma \circ
M_k(\phi(x)) \ \fa x \in M_+(\psi)} is called a
\textbf{\emph{Connes-Dixmier functional on $M_+(\psi)$}}.
\end{dfnS}

\begin{dfnS}
Let $\gamma \in SC_b^*([0,\iny))$ and $C$ be the Cesaro operator
of Section 3.1. We call a positive linear functional on
$C_b([0,\iny))$ of the form
$$
L_\gamma := \gamma \circ C
$$
a \textbf{\emph{Cesaro-Banach limit on $C_b([0,\iny))$}}. Let
$CBL(\RR_+)$ denote \textbf{\emph{the set of Cesaro-Banach limits
on $C_b([0,\iny))$}}.
\end{dfnS}

\REMS{The proof of Theorem 3.3 demonstrates that a Cesaro-Banach
limit $L_\gamma$ has the property of translation invariance and
$L_\gamma(1) = 1$.  Hence $L_\gamma \in BL(\RR_+)$ for all $\gamma
\in SC_b^*([0,\iny))$ and the set of Cesaro-Banach limits is a
proper subset of the set $BL(\RR_+)$,
$$CBL(\RR_+) \subset BL(\RR_+).$$
}

\noindent Let $k \in \mathcal{K}$. Define the continuous bounded
function \display{g_k(t) := g(k(t))} for any $t > 0$ and $g \in
C_b([0,\iny))$. Clearly, $g\to g_k$ is a $*$-automorphism of
$C_b([0,\iny))$.

Let $\gamma \in SC_b^*([0,\iny))$.  Then the functional $\gamma_k$
on $C_b([0,\iny))$ defined by \display{\gamma_k(g) := \gamma(g_k)
\ \fa g \in C_b([0,\iny))} has the properties $\gamma_k(1)=1$ and
$\gamma_k(f) = 0$ for all $f \in C_0([0,\iny))$. Hence $\gamma_k$
is an element of the set $SC_b^*([0,\iny))$.

\begin{propS} Let $k \in\mathcal{K}$. Then
$$
\gamma \circ M_k(g) = \gamma_{k} \circ C(g_{k^{-1}})
$$
for all $g \in C_b([0,\iny))$.
\par
\preskip
\proof{Using the substitution $s = k^{-1}(t)$,
$$
M_k(g)(\lambda) = \frac{1}{k(\lambda)}\int_0^\lambda g(s) dk(s)
=\frac{1}{k(\lambda)}\int_0^{k(\lambda)} g(k^{-1}(t)) dt
= C(g_{k^{-1}})(k(\lambda)).
$$
Hence $\gamma(M_k(g)) = \gamma_k(C(g_{k^{-1}}))$.
}
\end{propS}

\vspace*{0.1cm}

\begin{thmS}
Let $\psi \in \Omega_\iny$ and $k \in \mathcal{K}$.
\begin{prop2list}{10}{2}{4}
\item Let $\gamma \in SC_b^*([0,\iny))$. Then there
exists a Cesaro-Banach limit $L := L_{\gamma_k} \in CBL(\RR_+)$
such that
$$
\tau_{\gamma,k}(x) = f_{L,k^{-1}}(x) ,\quad \forall x \in
M_+(\psi).
$$
\item Let $L \in CBL(\RR_+)$.
Then there exists an element
$\gamma \in SC_b^*([0,\iny))$ such that
$$
f_{L,k^{-1}}(x) = \tau_{\gamma,k}(x),\quad \forall
x \in M_+(\psi).
$$
\end{prop2list}
\par
\preskip
\proof{Immediate from Proposition 5.5.}
\end{thmS}

\noindent The result implies the following important
identification. The set of Connes-Dixmier functionals arising from
the function $k \in \mathcal{K}$ is the set
$$
\inset{\tau_{\gamma,k}}{\gamma \in SC_b^*([0,\iny))} =
\inset{f_{L,k^{-1}}}{ L \in CBL(\RR_+)}. \eqno{(*)}
$$

\subsection{Subsets of Banach Limits and the Cesaro Limit Property}

Let $\kappa \in\mathcal{K}$. The identification $(*)$ above introduces
to our analysis the set of functionals
$$
\inset{ f_{L,\kappa} }{ L \in \Lambda }
$$
where $\Lambda$ is a subset of $BL(\RR_+)$. We consider, in this
section, a sufficient condition on a subset $\Lambda \subset
BL(\RR_+)$ such that the statements of Theorem 2.8 and Corollary
3.9 can be extended to the set of functionals $ \inset{
f_{L,\kappa} }{ L \in \Lambda } $.

\begin{dfnS}
Let $\Lambda \subset BL(\RR_+)$.  We say $\Lambda$ has the Cesaro
limit property if, for each $g \in C_b([0,\iny))$, \display{ \{ a
, b \} \subset \inset{L(g)}{L \in \Lambda} } where $
a=\liminf_{t\to\iny}C(g)(t)$ and $b=\limsup_{t\to\iny}C(g)(t).$
\end{dfnS}

\medskip \noindent Let $\kappa \in \mathcal{K}$ and $\Lambda \subset BL(\RR_+)$.
Define a seminorm on $M(\psi)$ by setting for $x \in M(\psi)$
$$
\nm{x}_{\kappa,\Lambda} := \sup \inset{f_{L,\kappa}(|x|)}{L \in
\Lambda}.
$$

\begin{thmS}
Let $\psi \in \Omega_\iny$ and
$\kappa \in \mathrm{D}(\psi)$.
Let $\Lambda
\subset BL(\RR_+)$ have the Cesaro limit property.
Then there exists $0 < c < 1$ such that
$$
c \rho_1(x) \leq \nm{x}_{\kappa,\Lambda} \leq
\rho_1(x),\quad \forall x \in M(\psi).
$$
\par
\preskip
\proof{As
$$
\rho_1(x) = \limsup_{t \to \iny} \phi(x)(t)
= \limsup_{t \to \iny} \phi_\kappa(x)(t).
$$
there exists a sequence of positive numbers $\{ t_k \}_{k=1}^\infty$ with $t_k \to \iny$ as $k \to \iny$
such that
$$
\lim_{k \to \iny} \frac{\sigma(t_k)}{t_k} = \rho_1(x) \eqno{(5.1)}
$$
where $\sigma(t) = t \phi_\kappa(x)(t) \ , \ t > 0$.  We write
$\sigma(t) = \frac{t}{f(t)} \int_0^{\kappa(t)} x^*(s)ds$ where
$f(t) = \psi \circ \kappa(t)$. Then $\sigma'(t) =
(1-t\frac{f'(t)}{f(t)}) \phi_\kappa(x)(t) +  \frac{t}{f(t)}
\kappa'(t) x^*(\kappa(t))$. The hypothesis on $\kappa$ implies
$t\frac{f'(t)}{f(t)} < H$ for some $H \geq 1$ and all $t > 1$.
Hence
$$\sigma'(t) > (1-H) \phi_\kappa(x)(t) $$
Let $s \in [t_k,et_k]$.  Then
$$
\sigma(s) - \sigma(t_k) = \int_{t_k}^s \sigma'(t)dt
> (1-H) \int_{t_k}^s \phi_\kappa(x)(t)dt
\geq (1-H) \int_{t_k}^{et_k} \phi_\kappa(x)(t)dt
$$
since $1-H \leq 0$, and
\begin{eqnarray*}
\frac{1}{et_k} \int_{t_k}^{et_k} \frac{\sigma(s)}{s}ds
& > & \frac{1}{et_k} \Big( \sigma(t_k) + (1-H)\int_{t_k}^{et_k} \phi_\kappa(x)(t)dt \Big)
\int_{t_k}^{et_k} \frac{ds}{s} \\
& = & \frac{\sigma(t_k)}{et_k} + (1-H)\frac{1}{et_k}\int_{t_k}^{et_k} \phi_\kappa(x)(t)dt \\
& \geq & \frac{1}{e} \frac{\sigma(t_k)}{t_k} + (1-H)C(\phi_\kappa(x))(et_k)
\end{eqnarray*}
as $\int_{t_k}^{et_k} \frac{ds}{s} = \log\frac{et_k}{t_k} = \log e = 1$.
Combining the previous inequality with
$$
C(\phi_\kappa(x))({et_k})
= \frac{1}{et_k} \int_0^{et_k}
\frac{\sigma(s)}{s} ds
\geq \frac{1}{et_k} \int_{t_k}^{et_k}
\frac{\sigma(s)}{s} ds
$$
yields
$$
C(\phi_\kappa(x))({et_k}) > \frac{1}{e} \frac{\sigma(t_k)}{t_k} + (1-H)C(\phi_\kappa(x))(et_k).
$$
Hence, after rearrangement,
$$
H \, C(\phi_\kappa(x))({et_k}) > \frac{1}{e} \frac{\sigma(t_k)}{t_k}
$$
and
$$
\limsup_{t \to \iny} C(\phi_\kappa(x))(t) \geq \limsup_{k \to
\iny} C(\phi_\kappa(x))({et_k}) > \frac{1}{He} \limsup_{k \to
\iny} \frac{\sigma(t_k)}{t_k} \stackrel{(5.1)}{=} c \rho_1(x)
$$
where $c = (eH)^{-1}$.  By the Cesaro limit property, $\exists L \in \Lambda$ such that
$L(\phi_\kappa(x)) = \limsup_{t \to \iny} C(\phi_\kappa(x))(t)$.
Hence $f_{L,\kappa}(x) > c \rho_1(x)$ for some $L \in
\Lambda$. The reverse inequality $f_{L,\kappa}(x) \leq \rho_1(x)$
for all $L \in \Lambda$ is obvious. }
\end{thmS}

\bigskip \noindent We now extend the notion of measurability and Definition 3.5.
Let $\Lambda \subset BL(\RR_+)$.  Define the set
\display{\mathcal{F}_\Lambda = \inset{g \in C_b([0,\iny))}{L_1(g)
= L_2(g) \fa L_1,L_2 \in \Lambda}. } Let $g \in
\mathcal{F}_\Lambda$.  We denote the value $A = L(g) \fa L \in \Lambda$ by
$$
\text{F$_\Lambda$-}\lim_{t \to \iny} g(t) = A.
$$

\begin{dfnS}
Let $\psi \in \Omega_\iny$ and $\kappa\in \mathcal{K}$. Let $x \in
M_+(\psi)$.  We say $x$ is
$\mathcal{F}_{\Lambda,\kappa}$-measurable if $\phi_\kappa(x) \in
\mathcal{F}_\Lambda$.
\end{dfnS}

\begin{thmS}
Let $\psi \in \Omega_\iny$ and
$\kappa \in \mathrm{D}(\psi)$.  Let $\Lambda \subset BL(\RR_+)$
have
the Cesaro limit property.
Then the following statements are equivalent
\begin{prop2list}{16}{2}{3}
\item $x$ is $\mathcal{C}_\kappa$-measurable, \item $x$ is
$\mathcal{F}_\kappa$-measurable, \item $x$ is
$\mathcal{F}_{\Lambda,\kappa}$-measurable, \item $x$ is
$\mathcal{S}$-measurable.
\end{prop2list}
\par
\preskip \proof{(iv) $\implies$ (ii) $\implies$ (iii) is
immediate. (iii) $\implies$ (i) is immediate from the Cesaro limit
property. (i) $\implies$ (ii) $\implies$ (iv) by Corollary 3.9. }
\end{thmS}

\subsection{Results on Connes-Dixmier Functionals}

We now concentrate
on the subset $CBL(\RR_+) \subset BL(\RR_+)$.

\begin{propS}
The set of Cesaro-Banach limits $CBL(\RR_+)$
has the Cesaro limit property.
\par
\preskip
\proof{The proof of Theorem 3.3.}
\end{propS}

\noindent With the identification $(*)$ of Section 5.1,
the results of Section 5.2 can be applied to the set
of Connes-Dixmier functionals as follows.

\begin{thmS}
Let $\psi \in \Omega_\iny$ and $k^{-1} \in \mathrm{D}(\psi)$. Then
\begin{prop2list}{10}{2}{4}
\item[\text{A.}] for each $x \in M(\psi)$,
$$
\rho_1(x) \simeq \sup \inset{\tau_{\gamma,k}(|x|)}{\gamma \in
SC_b^*([0,\iny))};
$$
\item[\text{B.}]
the following statements are equivalent
\begin{prop2list}{10}{2}{4}
\item $x$ is $\mathcal{C}_{k^{-1}}$-measurable, \item
$\tau_{\gamma,k}(x)$ is independent of $\gamma \in
SC_b^*([0,\iny))$, \item $$\tau_{\gamma,k}(x) = \lim_{t \to \iny}
\frac{1}{\psi(t)} \int_0^t x^*(s)ds \ \ \fa \gamma \in
SC_b^*([0,\iny)).$$
\end{prop2list}
\end{prop2list}
\par
\preskip
\proof{
Proposition 5.11, Theorem 5.8 and Theorem 5.10.
}
\end{thmS}

\REMS{We recall for the reader the particular
case of Connes' construction in \cite{CN}.
The pair of functions used in \cite{CN} is
$(\psi(t), k(t)) = (\log(1+t),\log(1+t))$.  It is
trivial to check
$k^{-1} \in \mathrm{D_{exp}}(\psi)$ and hence satisfies
the hypothesis of Theorem 5.12.

We note that the claim contained in Theorem 5.12 B.
(i) $\Leftrightarrow$ (ii)
generalises to arbitrary Marcinkiewicz spaces the assertion
proved by A. Connes for the choices $(\psi(t),k(t)) = (\log(1+t),\log(1+t))$
\cite[Proposition IV.2.$\beta$.6]{CN}.
The claim in Theorem 5.12 B. (ii) $\Leftrightarrow$ (iii) is new even for
$\psi(t) = \log(1+t)$.
}

\section{Application to Non-Commutative Geometry}

We conclude the paper by reducing the results to the setting of
singular traces on semifinite von Neumann algebras \cite{DPSS},
which includes, as the type $I$ case, the setting for
non-commutative geometry \cite[VI.2]{CN}.

\medskip We introduce the notation of \cite{DPSS} Section 4.
Let $(\mathcal{N},\tau)$ be the pair of a semifinite von Neumann
algebra $\mathcal{N}$ with a faithful normal semifinite trace
$\tau$.  Let $\chi_E$ denote the characteristic function of a
measurable set $E \subset [0,\iny)$. Define the generalised
singular values of the operator $r \in \mathcal{N}$ with respect
to $\tau$ \cite{FK}, \display{\mu_t(r) = \inf \inset{s \geq
0}{\tau(\chi_{(s,\iny)}(|r|)) \leq t} .} The function $\mu_t(r) :
[0,\iny) \to [0,\iny)$ is a non-increasing and right continuous
function. Define the Marcinkiewicz space $M(\log(1+t))$ as the set
of Lebesgue measurable functions $x : [0,\iny) \to [0,\iny)$
such that
$$
\nm{x}_{M(\log(1+t))} := \sup_{t > 0} \, \frac{1}{\log(1+t)} \int_0^t x^*(s) ds < \infty
$$
where $x^*$ is the decreasing rearrangement of $x$, see Section
1.2. Define $M_1(\log(1+t))$ as the closure of $L^1([0,\iny)) \cap
M(\log(1+t))$ in the norm $\nm{.}_{M(\log(1+t))}$. Define the
Marcinkiewicz (normed) operator ideal associated to the
Marcinkiewicz space $M(\log(1+t))$ by \display{
\begin{array}{l}
\mathcal{L}^{(1,\iny)}(\mathcal{N},\tau) := \inset{r \in \mathcal{N}}{ \mu_t(r) \in M(\log(1+t))} \\
\mathcal{L}_0^{(1,\iny)}(\mathcal{N},\tau) := \inset{r \in \mathcal{N}}{
\mu_t(r) \in M_1(\log(1+t))}
\end{array}
}
with norm
$$
\nm{r}_{(1,\iny)} := \nm{\mu_t(r)}_{M(\log(1+t))} \ \text{ for } \
r \in \mathcal{L}^{(1,\iny)}(\mathcal{N},\tau).
$$
We note the separable ideal $\mathcal{L}_0^{(1,\iny)}(\mathcal{N},\tau)$ is
the closure in the norm $\nm{.}_{(1,\iny)}$ of the ideal $
\mathcal{L}^1(\mathcal{N},\tau)$ of all $\tau$-integrable elements from
$\mathcal{N}$.

\bigskip We recall that $\mathcal{K}$ is the set of strictly increasing, invertible,
differentiable and unbounded functions mapping $[0,\iny) \to [0,\iny)$ and $BL(\RR_+)$
is the set of translation invariant positive linear functionals on $C_b([0,\iny))$, see Section 1.1 and Remark 1.9.
Let $k \in \mathcal{K}$ and $L \in BL(\RR_+)$.  Define a
functional on $\mathcal{L}^{(1,\iny)}(\mathcal{N},\tau)$ by
\display{ F_{L,k}(r) := L \Big( \frac{1}{\log(1 + k^{-1}(t))} \int_0^{k^{-1}(t)} \mu_s(r)ds \Big) }
for all positive
elements $r \in \mathcal{L}^{(1,\iny)}(\mathcal{N},\tau)$. We
extend $F_{L,k}$ to the positive part of $\mathcal{N}$ by setting
$F_{L,k}(r)=\infty$ for all positive elements $r\in
\mathcal{N}\backslash \mathcal{L}^{(1,\iny)}(\mathcal{N},\tau)$. A
linear functional $F$ on the von Neumann algebra $\mathcal{N}$ is
called singular (with respect to the faithful normal semi-finite
trace $\tau$) if $F$ vanishes on $\mathcal{L}^1(\mathcal{N},\tau)\cap
\mathcal{N}$.  We recall the notation F-$\lim_{n \to \iny} a_n = A$,
introduced by G. Lorentz in \cite{GL}, denotes almost convergence of
a sequence $\{ a_n \}_{n \in \NN}$ to the value $A \in \RR$, see Definition 1.1.

\begin{thmS}[Trace Theorem]
\medskip \noindent Let $(\mathcal{N},\tau)$ be a semifinite von Neumann factor $\mathcal{N}$
with faithful normal semifinite trace $\tau$. Then $F_{L,k}$ is a
singular trace on $\mathcal{N}$ if and only if $k^{-1}$ satisfies
\begin{prop2list}{10}{2}{3}
\item \text{F-}$\lim_{n \to \iny} \frac{\log(k^{-1}(n))}{\log(k^{-1}(n+1))} = 1$,
\item $\exists \, C > 0$ such that $k^{-1}(t + C) > 2k^{-1}(t)$ for all $t > 0$.
\end{prop2list}
\par \preskip \proof{Let $\alpha := \log(1+t)$. By construction $F_{L,k}(r) := f_{L,k^{-1}}(\mu_t(r))$,
where $f_{L,k^{-1}}$ is given in Definition 1.7,
and condition (i) and (ii) are equivalent to
$k^{-1} \in \mathrm{R_{exp}}(\alpha)$, see Definitions 2.1,2.5.
The
functional $f_{L,k^{-1}}$ is a positive homogeneous functional on
\mbox{$M_+(\alpha)$} satisfying (i) and (ii) of Definition 1.5.
We claim $F_{L,k}$ is a singular trace on $\mathcal{N}$ if and
only if $f_{L,k^{-1}} \in M_+(\alpha)^*_{\mathrm{sym},\iny}$.

$(\Rightarrow)$ The functional $F_{L,k}$ is a singular trace on
$\mathcal{N}$.  Hence it is additive by hypothesis.  Let
$\mathcal{N}$ be a type II (respectively, I) factor. Then by
Theorem 4.4 (respectively, Theorem 4.5) from \cite{DPSS} the
functional $f_{L,k^{-1}}$ on $M(\alpha)$ is additive. Hence
$f_{L,k^{-1}} \in M_+(\alpha)^*_{\mathrm{sym},\iny}$.

$(\Leftarrow)$ The functional $f_{L,k^{-1}}$ on $M(\alpha)$ is
symmetric by hypothesis.  Then by Theorem 4.2 \cite{DPSS} the
functional $F_{L,k}$ as a functional on
$\mathcal{L}^{(1,\iny)}(\mathcal{N},\tau)$ is additive. Let $r \in
\mathcal{L}^{(1,\iny)}(\mathcal{N},\tau)$ be positive and $u \in
\mathcal{N}$ be unitary. Then
$$F_{L,k}(uru^*) = f_{L,k^{-1}}(\mu_t(uru^*))
        = f_{L,k^{-1}}(\mu_t(r)) = F_{L,k}(r).
$$
Hence $F_{L,k}$ defines a trace.  The fact that
$F_{L,k}$ is a singular trace is immediate.

The result follows as $f_{L,k^{-1}} \in
M_+(\alpha)^*_{\mathrm{sym},\iny}$ if and only if $k^{-1} \in
\mathrm{R_{exp}}(\alpha)$ by Theorem 2.8.
}
\end{thmS}

\bigskip \noindent Define the dilation operator
\display{D_a(g)(b) = g(ab) \ \fa a,b \in (0,\iny), g \in
C_b([0,\iny)).} An element $\omega \in C_b([0,\iny))^*$ is called
dilation invariant if \display{\omega(D_{a}(g)) = \omega(g) \fa a
\in (0,\iny), g \in C_b([0,\iny)).} We recall $SC_b^*([0,\infty))$
denotes the set of all positive linear functionals $\gamma$ on
$C_b([0,\iny))$ such that $\gamma(1)=1$ and $\gamma(f) = 0$ for
all $f$ in $C_0([0,\iny))$, see Definition 5.1.  Define
$$
D(\RR_+) := \inset{\omega \in SC_b^*([0,\infty))}{  \ \omega
\text{ is dilation invariant}}.
$$
Let $\omega \in D(\RR_+)$. Define the functional $Tr_{\omega}$  on
$\mathcal{L}^{(1,\iny)}(\mathcal{N},\tau)$ by setting
\display{Tr_{\omega}(r) := \omega \Big( \frac{1}{\log(1+t)} \int_0^t \mu_s(r) ds \Big) }
for all positive $r \in \mathcal{L}^{(1,\iny)}(\mathcal{N},\tau)$. We shall refer to the
functional $Tr_{\omega}$ as a Dixmier trace.

\medskip \noindent Let $\alpha(t) = \log(1+t)$. Define
$$
M_\alpha(g)(\lambda) := \frac{1}{\log(1+\lambda)}\int_0^\lambda
g(s) d\log(1+t), \ \lambda>0$$ for $g \in C_b([0,\iny))$ as in
Section 5.1.  Let $\gamma \in SC_b^*([0,\iny))$. Define the
functional $tr_{\gamma}$ on
$\mathcal{L}^{(1,\iny)}(\mathcal{N},\tau)$ by setting
$$
tr_{\gamma}(r) := \gamma \circ M_\alpha
\Big( \frac{1}{\log(1+t)} \int_0^t \mu_s(r) ds \Big)
$$
for all positive $r \in \mathcal{L}^{(1,\iny)}(\mathcal{N},\tau)$.
We shall refer to the functional $tr_{\gamma}$ as
a Connes-Dixmier trace, after its introduction and use by A. Connes in [2].

\medskip \noindent After Cesaro and Hardy \cite[Section 1.3]{Hardy}, define the Cesaro mean by
$$
C(\lambda) = \frac{1}{\lambda} \int_0^\lambda g(s) ds, \ \lambda>0
$$
for $g \in C_b([0,\iny))$.  Let $\gamma \in SC_b^*([0,\iny))$. The
composition $\gamma \circ C$ is called a Cesaro-Banach limit, see
Definition 5.3, and the set of Cesaro-Banach limits, denoted
$CBL(\RR_+)$, is a proper subset of $BL(\RR_+)$.

\medskip  The identification of Dixmier and Connes-Dixmier traces corresponding
to the pair $(\mathcal{N},\tau)$ of a semifinite von Neumann algebra $\mathcal{N}$
and faithful normal semifinite trace is as follows.

\begin{thmS} Let $(\mathcal{N},\tau)$ be a semifinite von Neumann algebra $\mathcal{N}$
with faithful normal semifinite trace $\tau$. Then
\begin{eqnarray*}
\inset{Tr_\omega}{\omega \in D(\RR_+)}
& = & \inset{F_{L,\alpha}}{L \in BL(\RR_+)} \\[4pt]
\inset{ \, tr_\gamma \, }{\gamma \in SC_b^*([0,\iny))}
 & = & \inset{F_{L,\alpha}}{L \in CBL(\RR_+)}.
\end{eqnarray*}
where $\alpha(t) = \log(1+t)$.
\par
\preskip
\proof{
By construction $F_{L,k}(r) := f_{L,k^{-1}}(\mu_t(r))$,
where $f_{L,k^{-1}}$ is given in Definition 1.7.
By Remark 1.9, $f_{L,\log(1+t)^{-1}}=f_{L,e^t-1}\equiv
f_{L,\exp}$. Therefore, to prove the first equality, it is
sufficient to show that

(i) for a given $\omega \in D(\RR_+)$, there exists an $L\in
BL(\RR_+)$ such that
$$L(\phi(x)(e^t))=\omega(\phi(x)(t)),\quad 0\leq x\in \mathcal{L}^{(1,\iny)}(\mathcal{N},\tau);$$
where $\phi(x)(t) := (\log(1+t))^{-1} \int_0^t x(s)ds$.

(ii) for a given $L\in BL(\RR_+)$, there exists an $\omega \in
D(\RR_+)$ such that the equality above holds.

To establish (i), fix an $\omega \in D(\RR_+)$ and define
$L(g):=\omega(g_1(\log(t)))$, $g\in C_b([0,\infty))$, $t\ge 0$,
where we set $g_1(s):=0$ if $s<-1$, $g_1(s)=g(s)$ if $s\ge 0$,
$g_1(s)=(1+s)g(0),\ -1\le s<0$. Clearly, $g_1$ is continuous on
$\RR$.  We show that $L\in BL(\RR_+)$. It is evident that $L$ is a
positive linear functional on $C_b([0,\infty))$ which takes value
$1$ on $g(t)\equiv 1$ and vanishes on $C_0([0,\infty))$. Thus, it
remains to show that $L$ is translation invariant. Fix an $a>0$
and consider $L(T_a(g))=\omega((T_a(g))_1(\log(t))$. For all
sufficiently large $t>0$, the value $(T_a(g))_1(\log(t))$
coincides with $g(\log(t)+a)$. On the other hand, since $\omega
\in D(\RR_+)$, we have
$$\omega(g_1(\log(t))=\omega(D_{e^a}g_1(\log(t))=\omega(g_1(\log(t\cdot
e^a)))=\omega(g_1(\log(t)+a)).
$$
As the value $g_1(\log(t)+a)$ also coincides with $g(\log(t)+a)$
for all sufficiently large $t>0$, we conclude that
$L(T_a(g))=L(g)$.

To show (ii), fix an $L\in BL(\RR_+)$  and define
$\omega(g):=L(g(e^t))$, $g\in C_b([0,\infty))$. Again, it is clear
that $\omega$ is a positive linear functional on $C_b([0,\infty))$
which takes value $1$ on $g(t)\equiv 1$ and vanishes on
$C_0([0,\infty))$. To show that $\omega$ is dilation invariant,
fix an arbitrary $\lambda\ge 0$. The translation invariance of $L$
immediately yields that for every $r\in [0,\infty)$
$$
L(g(e^t))=L(T_r(g(e^t))=L(g(e^{t+r}))
$$
and so, setting $r:=-\min\{0,\log(\lambda)\}$, we obtain
$$
\omega(D_\lambda g)=L((D_\lambda g)(e^{t+r}))=L(g(\lambda
e^{t+r}))=L(g(e^{t+(r+\log(\lambda))}),\quad g\in C_b([0,\iny)).
$$
By construction $r+\log(\lambda)>0$ and again appealing to the
translation invariance of $L$, we conclude $ \omega(D_\lambda
g)=L(g(e^{t}))=\omega(g)$. This completes the proof of the first
equality.

The second equality follows from Theorem 5.6 as $tr_\gamma(r) =
\tau_{\gamma, \log(1+t)}(\mu_s(r))$ where $\tau_{\gamma,\log(1+t)}$
is a Connes-Dixmier functional, see Definition 5.2.
}
\end{thmS}

Theorem 6.2 completes the identification suggested by the results
of \cite{CPS2}.  It follows from Theorem 6.2:

\begin{thmS} Let $(\mathcal{N},\tau)$ be a semifinite von Neumann algebra $\mathcal{N}$
with faithful normal semifinite trace $\tau$. Then
\begin{prop2list}{10}{2}{3}
\item the functionals $Tr_\omega$ and $tr_\gamma$ on $\mathcal{L}^{(1,\iny)}(\mathcal{N},\tau)$
define singular traces on $\mathcal{N}$ for all
$\omega \in D(\RR_+)$ and $\gamma \in SC_b^*([0,\iny))$,
\item the set of Connes-Dixmier traces is a subset of the set of Dixmier traces.
\end{prop2list}
\par
\preskip \proof{The implication $(\Leftarrow)$ in the statement of
Theorem 6.1 does not require that $\mathcal{N}$ be a factor. Hence
(i) and (ii) follow from Theorem 6.2 and Example 4.3.}
\end{thmS}

Let $(\mathcal{N},\tau)$ be a semifinite von Neumann algebra $\mathcal{N}$
with faithful normal semifinite trace $\tau$.
The identification in Theorem 6.2 allows the results of previous sections
to be applied to the Marcinkiewicz operator ideal
$\mathcal{L}^{(1,\iny)}(\mathcal{N},\tau)$ as follows.

\begin{thmS}
Let $r \in \mathcal{L}^{(1,\iny)}(\mathcal{N},\tau)$.  Then
\display{\nm{r}_0 := \inf_{r' \in \mathcal{L}^{(1,\iny)}_0(\mathcal{N},\tau)}
\nm{r-r'}_{(1,\iny)}
        = \sup_{\omega \in D(\RR_+)} Tr_\omega(|r|)}
and
\display{\nm{r}_0 := \inf_{r' \in \mathcal{L}^{(1,\iny)}_0(\mathcal{N},\tau)}
\nm{r-r'}_{(1,\iny)}
        \simeq \sup_{\gamma \in SC_b^*([0,\iny))} tr_\gamma(|r|).}

\par
\preskip \proof{Remark 5.13, Theorem 5.12 A., Theorem 6.2 and
Theorem 2.8 since $\nm{r}_0 \equiv \rho_1(\mu_t(r))$. }
\end{thmS}

\bigskip \noindent In the
following definition (iii) follows A. Connes (see
\cite[IV.2.$\beta$, Proposition 6, Definition 7]{CN}).

\begin{dfnS}
Let $r \in \mathcal{L}^{(1,\iny)}(\mathcal{N},\tau)$ be positive.  Then we
say $r$ is
\begin{prop2list}{10}{2}{4}
\item $\mathcal{M}$-measurable if
$$\lim_{\lambda \to \iny} M_\alpha \Big(
\frac{1}{\log(1+t)}
\int_0^t \mu_s(r)ds \Big)(\lambda) = A
$$
for some $A \geq 0$, \item $\mathcal{F}$-measurable if
$Tr_\omega(r)$ is independent of $\omega \in D(\RR_+)$, \item
measurable if $tr_\gamma(r)$ is independent of $\gamma \in
SC_b^*([0,\iny))$, \item Tauberian if
$$\lim_{t \to \iny} \frac{1}{\log(1+t)}
\int_0^t \mu_s(r)ds = A$$
for some $A \geq 0$.
\end{prop2list}
\end{dfnS}

\begin{thmS}
Let $r \in \mathcal{L}^{(1,\iny)}(\mathcal{N},\tau)$ be positive. Then the
following statements are equivalent
\begin{prop2list}{10}{2}{4}
\item $r$ is $\mathcal{M}$-measurable,
\item $r$ is $\mathcal{F}$-measurable,
\item $r$ is measurable,
\item $r$ is Tauberian.
\end{prop2list}
\par
\preskip \proof{From the proof of Proposition 5.5 $\mu_t(r)$ is
$\mathcal{M}$-measurable if and only if $\mu_t(r)$ is
$\mathcal{C}_{\alpha^{-1}}$-measurable (see Definition 3.5). Hence
the result follows directly from Remark 5.13 and Theorem 5.10. }
\end{thmS}

\REMS{As mentioned in Remark 5.13, the equivalence of the
statements (i) and (iii) in Theorem 6.6 is a result stated and
proved by A. Connes \cite[IV.2.$\beta$, Proposition 6]{CN} for the
special case $(\mathcal{N},\tau) = (B(H),Tr)$ where $H$ is a
separable Hilbert space and $Tr$ is the canonical trace.  That a
positive element $r \in \mathcal{L}^{(1,\iny)}(\mathcal{N},\tau)$
is measurable if and only if $r$ is Tauberian is a new result. }

Theorem 6.6 has the following corollary, which shall conclude
the paper, linking measurable operators and
results of \cite{CPS2}.
\begin{corS}
Let $r \in \mathcal{L}^{(1,\iny)}(\mathcal{N},\tau)$ be positive.  Define
$$
\zeta_r(s) = \tau(r^s)
$$
for any $s \in \CC$ with Re$(s) > 1$.
Then the
following statements are equivalent
\begin{prop2list}{10}{2}{4}
\item $r$ is measurable,
\item $$
Tr_{\omega}(r) = \lim_{s \to 1^+} (s-1)\zeta_r(s)
= 2 \Gamma( \textstyle{\frac{1}{2}})^{-1} \lim_{\epsilon \to 0^+} \epsilon \tau(e^{-{(\epsilon r)}^{-2}})
$$
for all $\omega \in D(\RR_+)$.
\end{prop2list}
\par
\preskip \proof{Theorem 6.6, Corollary 3.7 \cite{CPS2} and
Proposition 4.2 \cite{CPS2}. }
\end{corS}

\noindent S.L. and F.S.: School of
Informatics and Engineering, Flinders University of South
Australia, Bedford Park, 5042, Australia

\noindent e-mail:
steven.lord@alumni.adelaide.edu.au, sukochev@infoeng.flinders.edu.au

\smallskip
\noindent A.S.: Department of Mathematics, Voronezh State
University of Architecture and Construction, 20-letiya Oktyabrya
84, Voronezh 394006, Russia

\noindent e-mail: sed@vmail.ru

\end{document}